\documentclass[11pt]{article}
\usepackage[dvips]{graphics}
\title{\Large\bf An edge crack and a crack close to the vertex of a wedge}
\author{\bf Y.A.\ Antipov\\ 
Department of Mathematics, Louisiana State University\\
Baton Rouge LA 70803, U.S.A.\\
yantipov@lsu.edu}

\date{}

\newcommand{\sgn}{\mathop{\rm sgn}\nolimits}

\newcommand{\I}{\mathop{\rm Im}\nolimits}
\newcommand{\R}{\mathop{\rm Re}\nolimits}
\newcommand{\const}{\mbox{const}}

\newcommand{\beqa}{\begin{eqnarray}}
\newcommand{\eeqa}[1]{\label{#1}\end{eqnarray}}
\newcommand{\bequ}{\begin{equation}}
\newcommand{\eequ}[1]{\label{#1}\end{equation}}

\newcommand{\Md}{\partial}

\newcommand{\Ga}{\alpha}
\newcommand{\Gb}{\beta}
\newcommand{\Gd}{\delta}
\newcommand{\Ge}{\epsilon}
\newcommand{\Gve}{\varepsilon}

\newcommand{\Gg}{\gamma}
\newcommand{\Gc}{\chi}

\newcommand{\Gl}{\lambda}
\newcommand{\Gn}{\eta}
\newcommand{\Gm}{\mu}

\newcommand{\Gt}{\theta}

\newcommand{\Gr}{\rho}

\newcommand{\Gs}{\sigma}

\newcommand{\Gj}{\tau}

\newcommand{\Go}{\omega}

\newcommand{\GD}{\Delta}

\newcommand{\GG}{\Gamma}
\newcommand{\GL}{\Lambda}

\newcommand{\GS}{\Sigma}

\newcommand{\GO}{\Omega}

\newcommand{\CP}{{\cal P}}

\newcommand{\beq}{\begin{equation}}
\newcommand{\eeq}{\end{equation}}
\newcommand{\barr}{\begin{eqnarray}}
\newcommand{\earr}{\end{eqnarray}}
\newcommand{\beqn}{\begin{equation*}}
\newcommand{\eeqn}{\end{equation*}}
\newcommand{\barrn}{\begin{eqnarray*}}
\newcommand{\earrn}{\end{eqnarray*}}

\newcommand{\fr}{\frac}

\begin{document}
\maketitle

\begin{abstract}

Two model problems of an elastic wedge with an internal and edge crack  are analyzed.  The problem of an internal crack
reduces
to an order-4 vector Riemann-Hilbert problem whose matrix kernel entries are meromorphic functions and have exponential factors.
When the internal crack is located along one of the wedge sides, an efficient method of solution is proposed. It requires a factorization
of the order-2 matrix coefficient associated with the corresponding problem of an edge crack and the solution of an infinite system 
of linear algebraic system with an exponential  rate of convergence of an approximate solution to the exact one. 
The order-2 Khrapkov's factorization
is modified by splitting the matrix kernel into a scalar dominant function and a ``regular" matrix whose factorization is more convenient for numerical purposes.
  Expressions for the stress intensity coefficients and the potential energy released when the crack advances are derived. Asymptotic relations for the  stress intensity coefficients and the potential energy when
one of the crack tips is close the wedge vertex 
are obtained.

\end{abstract}

\setcounter{equation}{0}

\section{Introduction}

The study of the effects of a crack in a half-plane  and a wedge has a long history. The vast majority of the investigations concern
the model problems of an edge crack  located either in a half-plane and orthogonal to its boundary or in a wedge and lying on it bisector.  Both
analytical and numerical solutions for these particular cases are available in the literature. 
First  solutions and expressions for the stress intensity factor for the model problem of an edge slit 
in a half-plane orthogonal to its boundary were recovered   \cite{koi1},  \cite{irw}, \cite{bue}, \cite{bow} by using different approximate techniques.
  Koiter \cite{koi1} reduced the model to a Wiener-Hopf problem and solved it approximately by factorizing not the kernel 
but a function that approximates it. Later \cite{koi2} he recovered an exact formula  for the stress intensity factor. The first exact solution to this model problem
of an edge crack was obtained by Wigglesworth \cite{wig} by reducing the problem to a  Wiener-Hopf equation different  from  \cite{koi1} and
factorizing the kernel exactly by means of  infinite products. Exact representations of the stress and displacement fields were found \cite{dor} by applying the method \cite{koi1}
and factorizing the kernel in terms of generalized factorial functions. 

The problem of an edge crack along the bisector of a wedge is also solvable by quadratures. Bantsuri  \cite{ban} used a conformal map 
and the Muskhelishvili method to transform the wedge with
an edge slit along the bisector into a strip cut along the negative semi-axis, derived a scalar Riemann-Hilbert problem and solved it by means of the Cauchy integrals. Khrapkov \cite{khr1}
constructed an exact solution to the same problem by applying the Mellin transform to the governing  boundary value problem of the model and solving the associtaed Riemann-Hilbert problem 
different from \cite{ban}. The  problem of an edge symmetric slit in a wedge was solved \cite{sri} by perturbing the exact solution  \cite{wig} for a half-plane, reducing the problem to a Fredholm  integral equation and solving it approximately. An asymptotic solution for large 
 $\mu=1/log\fr{b}{a}$, that is when the internal symmetric
 crack $\{a<r<b, \Gt=0\}$ is far from the vertex $r=0$
of a wedge, was found in \cite{sme}.

If the wedge sides are free of traction and the edge crack is not on the wedge bisector, then the problem reduces to an order-2 vector Riemann-Hilbert  problem.  In general, 
for an arbitrary wedge angle, the corresponding matrix of the problem cannot be factorized by quadratures by methods available in the literature, and the exact solution 
cannot be recovered. However, in three particular nonsymmetric cases \cite{khr2}, exact solutions may be determined. These cases are:
(i) the wedge consists of a half-plane $\{-\pi<\Gt<0\}$ and a wedge $\{0<\Gt<\Ga\}$, while the edge crack is located in the segment $\{0<r<b, \Gt=0\}$, 
(ii) the wedge is a half-plane, while the edge crack forms with the boundary an arbitrary angle $\Ga\in (0,\pi)$,
and (iii) the wedge
is the whole plane cut along the ray $\{0<r<\infty, \Gt=\pi\}$ and a finite segment $\{0<r<b, \Gt=\Ga\in(0,\pi)\}$. Khrapkov \cite{khr2}
proposed a powerful method of efficient factorization of the corresponding matrix coefficient and derived closed-form solutions to these three model problems including expressions for
the weight matrices for the stress intensity factors by quadratures and in a series form. For computations, the first two terms of the expansions were recovered, tabulated and proposed to be used for computations. 

The main goal of this paper to study (1) the model problem of an internal crack $\{a<r<b, \Gt=0\}$ in a wedge $\{0<r<\infty, -\Ga_1<\Gt<\Ga_2\}$ when the crack is close to the vertex of the 
wedge and (2) the limit case $a=0$ when it becomes an edge crack. In the former case the problem reduces to an order-4 vector Riemann-Hilbert problem whose entries 
are meromorphic functions and some of them have exponential factors.
For its solution we advance the method proposed in \cite{ant1} for an order-2 vector Riemann-Hilbert problem associated with the classical problem of an annulus-shaped stamp and later applied to different models  \cite{ant2}, \cite{ant3}, \cite{ant4},  \cite{ant5}
arising in fracture, contact,  and penetration mechanics. An iterative method for factorization of triangular matrices with exponential factors was proposed in \cite{kis}.

In Section 2 we state the problem of an internal crack in a wedge. We show that it is equivalent to an oder-4 vector Riemann-Hilbert problem in Section 3.
In Section 4 and 5 we deal with a wedge $\{-\pi<\Gt<\Ga\}$, $\Ga\in(0,\pi)$. Since a building block of the solution procedure requires a factorization of the $2\times 2$
matrix associated with the limit case $a=0$, in Section 4 we revisit the first Khrapkov model problem of an edge crack in a wedge along one of the sides, namely $\{0<r<b, \Gt=0\}$.
We also modify the original Khrapkov scheme by splitting the matrix coefficient into a dominant scalar responsible for the proper singular behavior of the stresses at the crack
tip $r=b$ and a ``regular" matrix whose factorization becomes more convenient for numerical purposes. In this section we derive the weight matrix for the stress intensity
factors and consider two cases of loading. The first case is standard in such fracture problems: the normal and tangential loads are constants. For the second case
we consider first the associated singular problem of elasticity for a wedge without a crack. It is known  \cite{che} that if the wedge angle is in the range $(\pi,2\pi)$ and the wedge
sides are free of traction, then there is an eigen-solution that decays at infinity as $r^{\mu-1}$, $\Gm\in(0,1)$. We study the model problem for a wedge and an edge crack
with wedge and crack faces being free of traction when at infinity the eigen-solution is applied. In Section 5 we assume that the crack is internal. We derive an analytical  solution
to the problem whose representation formulas posses some quadratures (Cauchy integrals) and exponentially convergent series whose coefficients are recovered by
solving infinite systems of linear algebraic equation of the second kind. The rate of convergence of an approximate solution to the exact solution is exponential.
We derive the stress intensity factors $K^\pm=(K_I^\pm, K^\pm_{II})$ at the crack tips $a$ and $b$, respectively, and show that the factors at the crack tip $a=0$ behave as
\beq
K^-\sim -\fr{1}{\sqrt{\Gd}\log \Gd} Q K^+_0, \quad \Gd=\fr{a}{b}\to 0,
\label{1.1}
\eeq
where $Q$ is a matrix that depends 
on the wedge angle and independent of the loads, and $K^+_0$ is the vector whose entries are the stress intensity factors at the tip $b$ when $a=0$. We also determine the potential energy released when the one of the crack tips advances and the second end is fixed.
In Section 6 we show that even in the general case, when a crack $\{a<r<b, \Gt=0\}$ is in a wedge
 $\{0<r<\infty, -\Ga_1<\Gt<\Ga_2\}$, formula (1.1) is valid. However, in contrast to the model analyzed in Section  5, when
the matrix $Q$ is reconstructed in an explicit form, in the general case the matrix $Q$ is expressed through implicitly defined factorization factor-matrices. What is common is that, as in Section 5,
it is independent of the loads.
In Section 7 we analyze the model of an internal crack orthogonal to the boundary of a half-plane. In this particular case, the formulas are simplified since
the governing Riemann-Hilbert problem is an order-2 problem and instead of  matrix factorization one needs to factorize a scalar function only.
We modify the solution we obtained earlier \cite{ant2} and in addition, determine the potential energy
for any values of $a$ including the case $a\to 0$.

\setcounter{equation}{0}

\section{Formulation}

We consider the plane problem of an elastic wedge $W=\{0<r<\infty, -\Ga_2<\Gt<\Ga_1\}$,
$0<\Ga_j<\pi$,
 whose boundary is free of traction
$$
\Gs_\Gt=0, \quad \tau_{r\Gt}=0, \quad \Gt=-\Ga_2, \quad 0<r<\infty,
$$
\beq
\Gs_\Gt=0, \quad \tau_{r\Gt}=0, \quad \Gt=\Ga_1, \quad 0<r<\infty,
\label{2.1}
\eeq
The wedge $W$ is perturbed by the presence of a plane crack $\{a<r<b, \Gt=0^\pm\}$
whose faces are subject to the loads
\beq
\Gs_\Gt=-p_\Gt(r), \quad \tau_{r\Gt}=-p_{r\Gt}(r), \quad \Gt=0^\pm, \quad a<r<b.
\label{2.2}
\eeq
The functions $p_r$ and $p_{r\Gt}$ are integrable on the interval $(a,b)$, and the stresses
$\Gs_\Gt$ and $\tau_{r\Gt}$ vanish at infinity
as $r^{-1}$.
We confine ourselves to considering the plane stress conditions. In the plane strain case,
the Young modulus $E$ and the Poisson ratio  $\nu$ are to be replaced by $E_1=E/(1-\nu^2)$
and $\nu_1=\nu/(1-\nu^2)$. 
The displacements are discontinuous through the crack faces, and we denote by
$$
\fr{1}{E}\Gc_1(r)=\fr{\Md u_\Gt(r,0^-)}{\Md r}-\fr{\Md u_\Gt(r,0^+)}{\Md r},
$$
\beq
\fr{1}{E}\Gc_2(r)=\fr{\Md u_r(r,0^-)}{\Md r}-\fr{\Md u_r(r,0^+)}{\Md r}
\label{2.3}
\eeq
the jumps of their tangential derivatives which are unknown on the segment $(a,b)$ and vanish
if $r<a$ or $r>b$. 
Due to their definition the functions $\Gc_1(r)$ and  $\Gc_2(r)$ have to meet the conditions
\beq
\int_a^b \Gc_1(r)dr=0, \quad \int_a^b \Gc_2(r)dr=0.
\label{2.3'}
\eeq 
In the wedge $W$, the stresses satisfy the equilibrium equations
$$
r\fr{\Md\Gs_r}{\Md r}+\fr{\Md\Gj_{r\Gt}}{\Md \Gt}+\Gs_r-\Gs_\Gt=0,
$$
\beq
\fr{\Md\Gs_\Gt}{\Md \Gt}+r\fr{\Md\Gj_{r\Gt}}{\Md r}+2\Gj_{r\Gt}=0,
\label{2.4}
\eeq
and the compatibility condition
\beq
\GD(\Gs_r+\Gs_\Gt)=0.
\label{2.5}
\eeq

\setcounter{equation}{0}

\section{Order-4 vector Riemann-Hilbert problem}

 It is convenient to  apply the Mellin transform 
\beq
(\Gs_\Gt^s,\tau_{r\Gt}^s,\Gn_\Gt^s,\Gn_r^s)(\Gt)=
\int_0^\infty\left(\Gs_\Gt,\tau_{r\Gt},E\fr{\Md u_\Gt}{\Md r}, E\fr{\Md u_r}{\Md r}\right)(r,\Gt)r^s dr
\label{2.6}
\eeq
to  equations (\ref{2.4}) and (\ref{2.5}) in the wedges $-\Ga_2<\Gt<0$ and $0<\Gt<\Ga_1$ 
 separately.
The resulting one-dimensional system of equations yields the following relations
between the vectors $\Gs^s=(\Gs_\Gt^s,\tau_{r\Gt}^s)$ and $\Gn^s=(\Gn_\Gt^s,\Gn_r^s)$
on the line $\Gt=0$ \cite{khr2}:
$$
\Gn^s(0^-)=\fr{H(s,-\Ga_2)}{d(s,-\Ga_2)}\Gs^s(0),\quad
\Gn^s(0^+)=\fr{H(s,\Ga_1)}{d(s,\Ga_1)}\Gs^s(0),
$$
\beq
 s\in\GO=\{\R s=\Go, -\infty<\I s<\infty\}, \quad \Go\in\left(-\fr{1}{2},0\right),
\label{2.7}
\eeq
where $H(s,\Gt)$ is the matrix given by
\beq
H(s,\Gt)=\left(\begin{array}{cc}
\sin 2\Gt s+s\sin 2\Gt & 2s(s-1)\sin^2\Gt-(1-\nu)d(s,\Gt)\\
-2s(s+1)\sin^2\Gt+(1-\nu)d(s,\Gt) & \sin 2\Gt s-s\sin 2\Gt\\
\end{array}
\right),
\label{2.8}
\eeq 
and
\beq
d(s,\Gt)=s^2\sin^2\Gt-\sin^2\Gt s.
\label{2.9}
\eeq
Aiming to derive a vector Riemann-Hilbert problem we subtract the expression for $\Gn^s(0^+)$
in (\ref{2.7}) from the one for $\Gn^s(0^-)$ and deduce
\beq
\Gc^c(0)=G(s)\Gs^s(0), \quad s\in\GO,
\label{2.10}
\eeq
where $\Gc^s(0)=\Gn^s(0^-)-\Gn^s(0^+)$, 
and the entries $G_{jk}(s)$ of the matrix $G(s)$ have the form
$$
G_{jj}(s)=-\fr{\sin 2\Ga_1 s-(-1)^js\sin 2\Ga_1}{d(s,\Ga_1)}-
\fr{\sin 2\Ga_2 s-(-1)^js\sin 2\Ga_2}{d(s,\Ga_2)},
$$
\beq
G_{j \, 3-j}(s)=2s[(-1)^{j-1} s-1]\left[-\fr{\sin^2\Ga_1}{d(s,\Ga_1)}+
\fr{\sin^2\Ga_2}{d(s,\Ga_2)}
\right], \quad j=1,2.
\label{2.12}
\eeq

To clarify the analytic properties of the vector-functions $\Gc^s(0)$ and $\Gs^s(0)$, we 
represent them as
$$
\Gc^s(0)=b^{s+1}\Gc^-(s)=b^{s+1}\Gd^{s+1}\Gc^+(s), \quad \Gd=\fr{a}{b}\in(0,1),
$$
\beq
\Gs^s(0)=b^{s+1}[\Gd^{s+1}\Gs^-(s)+g^-(s)+\Gs^+(s)],
\label{2.13}
\eeq
where $g^-(s)$ is the only one vector-function in (\ref{2.13}) that is known. It is 
expressed through the  Mellin transforms of  the loads applied to the crack faces 
\beq
g^-(s)(s)=\Gd^{s+1}g^+(s),\quad 
g^-(s)=\int_\Gd^1\Gs(br)r^s dr, \quad 
g^+(s)=\int_1^{1/\Gd}\Gs(ar)r^s dr.
\label{2.14}
\eeq
The integral representations
of the unknown vector-functions
in (\ref{2.13}) are
$$
\Gc^-(s)=\int_\Gd^1\Gc(br)r^s dr, \quad \Gc^+(s)=\int_1^{1/\Gd}\Gc(ar)r^s dr, 
$$
\beq
\Gs^-(s)=\int_0^1\Gs(ar)r^s dr, \quad \Gs^+(s)=\int_1^{\infty}\Gs(br)r^s dr.
\label{2.15}
\eeq
Here and further, $\Gs(r)=(\Gs_\Gt(r,0),\tau_{r\Gt}(r,0))$ and  $\Gc(r)=(\Gc_1(r,0),\Gc_2(r))$.
The integrals $g^\pm(s)$ and $\Gc^\pm(s)$ are entire vector-functions in the complex
$s$-plane and  holomorphic
everywhere in the half-planes $D^\pm$,  $D^+=\{\R s < \Go\}$ and $D^-=\{\R s > \Go\}$. 
The integrals
$\Gs^-(s)$ and $\Gs^+(s)$ are holomorphic vector-functions in the half-planes $D^-$
and $D^+$, respectively.
The relations  (\ref{2.10}) and (\ref{2.13}), when combined, constitute
the following vector Riemann-Hilbert problem with a block-triangular
matrix coefficient for two pairs of vector-functions:
$$
\Gc^+(s)=\Gd^{-s-1}\Gc^-(s),
$$
\beq
\Gs^+(s)=[G(s)]^{-1}\Gc^-(s)-\Gd^{s+1}\Gs^-(s)-g^-(s),\quad s\in\GO.
\label{2.16}
\eeq

\setcounter{equation}{0}

\section{Edge crack in a wedge along one of its sides}\label{edge}

We start with a particular case of the problem when first, 
 $a=0$ that is one of the crack tips coincides with the wedge vertex and second,
$\Ga_2=\pi$ and $\Ga_1=\Ga\in(0,\pi)$. This means that the crack continues one of the 
wedge sides.

Due to the first assumption the vector $\Gs^-(s)=0$. The order-4 vector Riemann-Hilbert is significantly  simplified and 
becomes an order-2 problem
\beq
\Gc^-(t)=G(t)[g^-(t)+\Gs^+(t)], \quad t\in \GO.
\label{3.1}
\eeq
Here, $\Gc^-(s)$ and $\Gs^+(s)$ are unknown vector-functions, 
while $g^-(s)$ is a prescribed vector. Their integral representations are 
\beq
\Gc^-(s)=\int_0^1\Gc(br)r^s dr, \quad \Gs^+(s)=\int_1^{\infty}\Gs(br)r^s dr,\quad
g^-(s)=\int_0^1\Gs(br)r^s dr,
\label{3.2}
\eeq
and when $0<r<1$,  $\Gs(br)=(-p_\Gt(br), -p_{r\Gt}(br))$ describes the load applied to the crack faces.
The second assumption simplifies the matrix $G(s)$, and  it has the structure  
\cite{khr2}
\beq
G(s)=\left(\begin{array}{cc}
b_1(s)+b_2(s)l & b_2(s) m_+(s)\\
b_2(s)m_-(s) & b_1(s)-b_2(s)l\\
\end{array}\right), 
\label{3.3}
\eeq
whose entries are 
\beq
b_1(s)=2\cot\pi s-\fr{\sin 2\Ga s}{d(s,\Ga)}, \quad
b_2(s)=-\fr{2s\sin\Ga}{d(s,\Ga)},
\quad l=\cos\Ga, \quad m_\pm(s)=(\pm s-1)\sin\Ga.
\label{3.4}
\eeq

\subsection{Exact solution to the order-2 vector Riemann-Hilbert problem}\label{fact}

The first step of the solution method is to split the matrix $G(s)$ into a dominant scalar 
that yields the singularity of the solution at infinity and a ``regular" matrix. We have 
\beq
G(s)=4\cot\pi s G_0(s),
\label{3.5}
\eeq
where
$$
G_0(s)=\left(\begin{array}{cc}
b^\circ_1(s)+b^\circ_2(s)l & b^\circ_2(s) m_+(s)\\
b^\circ_2(s)m_-(s) & b^\circ_1(s)-b^\circ_2(s)l\\
\end{array}\right), 
$$
\beq
b^\circ_1(s)=\fr{1}{2}\left(1-\fr{\tan\pi s\sin 2\Ga s}{2d(s,\Ga)}\right), \quad
b^\circ_2(s)=-\fr{s\tan\pi s \sin\Ga}{2d(s,\Ga)}.
\label{3.6}
\eeq 
We factorize the dominant function in terms of the $\GG$-functions 
\beq
\cot\pi s=\fr{K^+(s)}{K^-(s)}, \quad K^+(s)=-\fr{\GG(-s)}{\GG(1/2-s)}, \quad 
K^-(s)=\fr{\GG(1/2+s)}{\GG(1+s)},
\label{3.7}
\eeq
and apply the Khrapkov factorization \cite{khr2} for the matrix $G_0(s)$ 
\beq
G_0(t)=X^+(t)[X^-(t)]^{-1}=[X^-(t)]^{-1}X^+(t), \quad t\in \GO,
\label{3.8}
\eeq
where $X^\pm(t)=X(\Go\mp 0+i\tau)$, $-\infty<\tau<\infty$,
$$
X(s)=B(s)\left(\begin{array}{cc}
c_+(s) & s_+(s)\\
s_-(s)& c_-(s)\\
\end{array}
\right),
$$
$$
c_\pm(s)=\cosh[f^{1/2}(s)\Gb(s)]\pm\fr{l}{f^{1/2}(s)}\sinh[f^{1/2}(s)\Gb(s)],
$$
$$
s_\pm(s)=\fr{m_\pm(s)\sinh[f^{1/2}(s)\Gb(s)]}{f^{1/2}(s)}, \quad f(s)=l^2+m_+(s)m_-(s),
$$
\beq
B(s)=\exp\left\{\fr{1}{4\pi i}\int_\GO\fr{\log\GD(t)dt}{t-s}\right\},
\quad
\Gb(s)=\fr{1}{2\pi i}\int_\GO\fr{\Ge(t)dt}{f^{1/2}(t)(t-s)},
\label{3.9}
\eeq
the function $\GD(s)$ is the determinant of the matrix $G_0(s)$, and this function and the function
$\Ge(t)$ are expresses through the eigenvalues $\Gl_1$ and $\Gl_2$ of the matrix $G_0(t)$
\beq
\Gl_1(t)=b_1^\circ(t)+b_2^\circ(t)f^{1/2}(t),\quad \Gl_2(t)=b_1^\circ(t)-b_2^\circ(t)f^{1/2}(t),
\label{3.10}
\eeq
as follows: 
\beq
\GD(s)=\Gl_1(s)\Gl_2(s), \quad \Ge(t)=\fr12\log\fr{\Gl_1(t)}{\Gl_2(2)}.
\label{3.11}
\eeq
The function $ f^{1/2}(s)=(1-s^2\sin^2\Ga)^{1/2}$ is a single branch of the two-values function
$w^2=1-s^2\sin^2\Ga$. The branch is  holomorphic in the $s$-plane cut along the line joining the branch points
$-1/\sin\Ga$ and  $1/\sin\Ga$ and passing through the infinite point. The factorization formulas
in  (\ref{3.9}) are independent of the branch  choice. 

The meromorphic functions $b_1^\circ(s)$ and $b_2^\circ(s)$ are even. 
At infinity,  $s=\Go+i\tau$, these functions behave as
\beq
b_1^\circ(s)\sim 1-4s^2\sin^2\Ga e^{2i\Ga s\sgn \tau}, 
\quad 
b_2^\circ(s)\sim -2is\sgn \tau \sin\Ga  e^{2i\Ga s\sgn \tau}, \quad  \tau\to\pm\infty.
\label{3.12}
\eeq
They have a
removable singularity at the point $s=0$,  and 
\beq
\Gl_j(s)\sim\fr12\left(1+\fr{\pi}{\Ga+(-1)^j\sin\Ga}
\right)>0, \quad s\to 0.
\label{3.12'}
\eeq
 It is convenient to shift the contour $\GO$ to the imaginary axis 
in the integral 
representations in (\ref{3.9}) and transform the integrals to the form
\beq
B(s)=\exp\left\{-\fr{s}{2\pi}\int_0^\infty\fr{\log\GD(i\tau)d\tau}{\tau^2+s^2}\right\},
\quad
\Gb(s)=-\fr{s}{\pi}\int_0^\infty\fr{\Ge(i\tau)d\tau}{\sqrt{1+\tau^2\sin^2\Ga}(\tau^2+s^2)}.
\label{3.13}
\eeq
The functions  $\log\GD(i\tau)$ and $\Ge(i\tau) $are bounded at $\tau=0$. At infinity, they behave as 
\beq
\log\GD(i\tau)\sim 8\tau^2\sin^2\Ga e^{-2\Ga\tau}, \quad 
\Ge(i\tau)\sim 2\tau^2\sin^2\Ga e^{-2\Ga\tau}, \quad \tau\to\infty.
\label{3.14}
\eeq
For real $s$, the integrals (\ref{3.13}) are real and rapidly converge.

After factorization of the function $\cot\pi s$ and matrix $G_0(s)$ is completed
we have 
\beq
K^-(t)X^-(t)\Gc^-(t)=4K^+(t)X^+(t)[\Gs^+(t)+g^-(t)], \quad t\in\GO.
\label{3.15}
\eeq
The Sokhotski-Plemelj formulas, when applied to the integral
\beq
\Psi(s)=\fr{1}{2\pi i}\int_\GO \fr{K^+(t)X^+(t)g^-(t)dt}{t-s},
\label{3.16}
\eeq
represent the vector-function $K^+(t)X^+(t)g^-(t)$ in terms of the limit values of the vector-function $\Psi(s)$
from the left and the right, $\Psi^+(t)$
and  $\Psi^-(t)$, respectively, 
\beq
K^+(t)X^+(t)g^-(t)=\Psi^+(t)-\Psi^-(t).
\label{3.17}
\eeq
The continuity principle, analysis of the functions at infinity and the Liouville theorem bring us to 
exact representation formulas for the problem solution. We obtain
$$
\Gc^-(s)=-\fr{4}{K^-(s)}[X^-(s)]^{-1}\Psi^-(s),\quad s\in D^-,
$$
\beq
\Gs^+(s)=-\fr{1}{K^+(s)}[X^+(s)]^{-1}\Psi^+(s),\quad s\in D^+.
\label{3.18}
\eeq

\subsection{Weight matrix}

Introduce the stress intensity factors at the crack tip $r=b$ in the standard way
\beq
\Gs_\Gt(r,0)\sim\fr{K_{I}}{\sqrt{2\pi(r-b)}}, \quad 
\tau_{r\Gt}(r,0)\sim\fr{K_{II }}{\sqrt{2\pi(r-b)}}, \quad r\to b+0.
\label{3.19}
\eeq 
On applying the abelian theorem to the integral $\Gs^+(s)$ given by (\ref{3.2})
we obtain their asymptotics at infinity
\beq
\Gs^+(s)\sim \fr{(-s)^{-1/2}}{\sqrt{2b}}\left(\begin{array}{c}
K_{I}\\
K_{II}\\
\end{array}
\right), \quad s\to\infty, \quad s\in D^+.
\label{3.20}
\eeq
On the other side, the asymptotics of the vector-function $\Gs^+(s)$ is defined from the analysis
of the solution (\ref{3.18})
\beq
\Gs^+(s)\sim  (-s)^{-1/2} X_\infty^{-1}\Psi_0,
 \quad s\to\infty, \quad s\in D^+.
 \label{3.21}
 \eeq
Here, 
\beq
\Psi_0=\fr{1}{2\pi i}\int_\GO K^+(t)X^+(t)g^-(t)dt,
\label{3.22}
\eeq
$X(s)\sim X_\infty$ as  $s\to \infty$, and since 
\beq
B(s)\sim 1, \quad c_\pm(s)\sim \cos q, \quad s_\pm(s)\sim\mp \sin q, \quad s\to\infty,
\label{3.23}
\eeq
where
\beq
q=\fr{\sin\Ga}{\pi}\int_0^\infty \fr{\Ge(i\tau)d\tau}{\sqrt{1+\tau^2\sin^2\Ga}},
\label{3.28}
\eeq
we have
\beq
X_\infty^{-1}=
\left(\begin{array}{cc}
\cos q & \sin q\\
-\sin q & \cos q\\
\end{array}
\right).
\label{3.29}
\eeq
By comparing formulas (\ref{3.20}) and (\ref{3.21}) and employing the third formula in (\ref{3.2})
 it is straightforward to express the stress intensity factors through the weight matrix
\beq
W(r)=X_\infty^{-1}\fr{1}{\pi i}\sqrt{\fr{b}{2}}\int_\GO K^+(t)X^+(t)r^tdt.
\label{3.30}
\eeq
We have
\beq
\left(\begin{array}{c}
K_{I}\\
K_{II}\\
\end{array}
\right)=\int_0^1 W(r)\Gs(br)dr,
\label{3.31}
\eeq
where $\Gs(r)=(-p_\Gt(r), -p_{r\Gt}(r))$ is the loading applied to the crack faces.

\subsection{Constant loading}\label{const}

\begin{figure}[t]
\centerline{
\scalebox{0.6}{\includegraphics{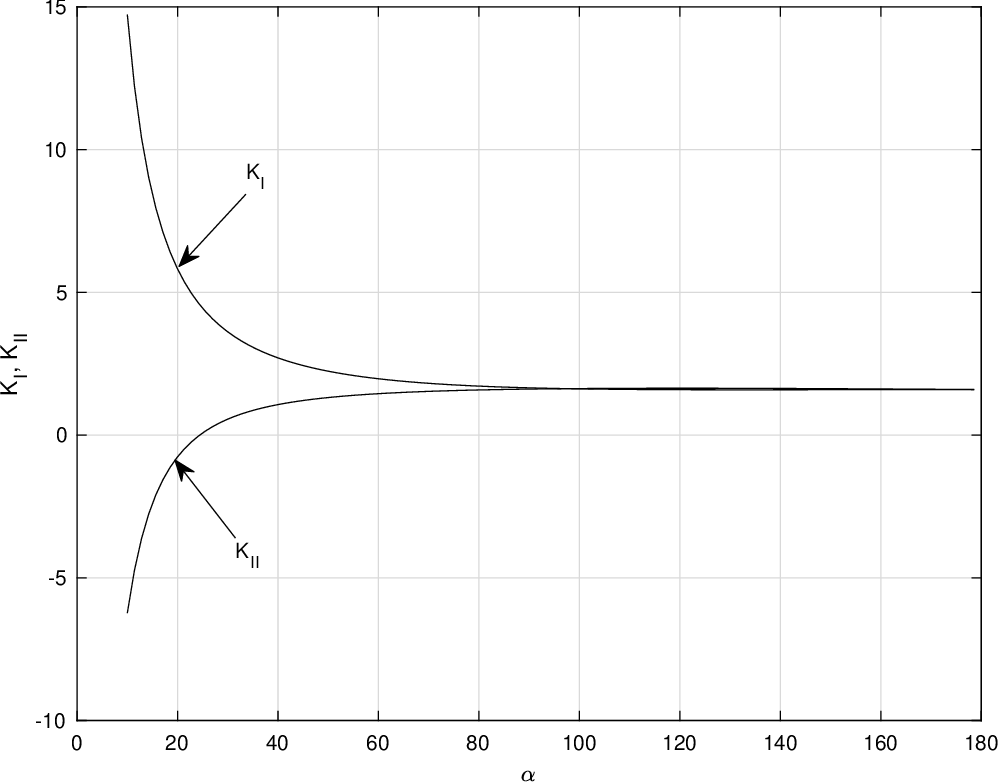}}}
\caption{
Stress intensity factors $K_I$ and $K_{II}$  versus the angle $\Ga$  measured in degrees when $b=1$ and $P_1=P_2=1$.}
\label{fig1}
\end{figure}

\begin{figure}[t]
\centerline{
\scalebox{0.6}{\includegraphics{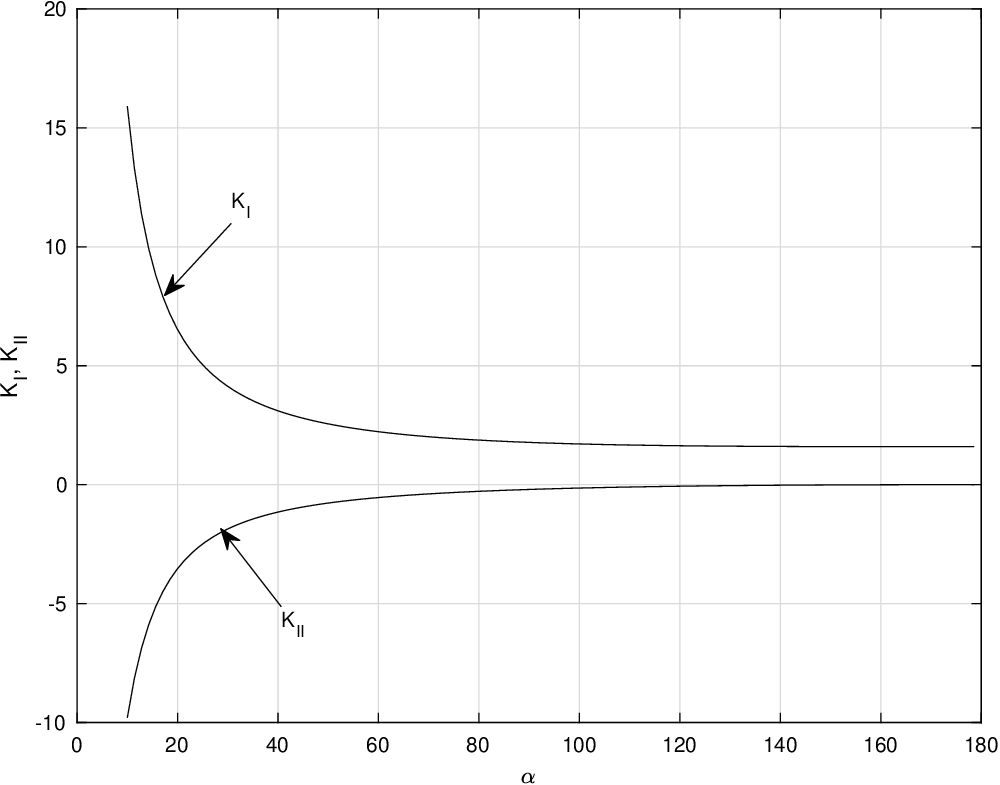}}}
\caption{
Stress intensity factors $K_I$ and $K_{II}$  versus the angle $\Ga$  measured in degrees  when $b=1$ and the loads are constant, $P_1=1$ and $P_2=0$.}
\label{fig2}
\end{figure} 

Consider now the case when both functions $p_\Gt(r)$ and $p_{r\Gt}(r))$ are constants,
$p_\Gt(r)=P_1$ and  $p_{r\Gt}(r)=P_2$. Then the vector-function $g^-(s)$ becomes
\beq
g^-(s)=-\fr{1}{s+1}P, \quad P= \left(\begin{array}{c}
P_1\\
P_2\\
\end{array}
\right).
\label{3.32}
\eeq
In this case the integral (\ref{3.16}) can be bypassed, and the solution 
of the vector Riemann-Hilbert problem (\ref{3.1}) has a simpler form
$$
\Gc^-(s)=\fr{8}{\sqrt{\pi}(s+1)K^-(s)}[X^-(s)]^{-1}X^+(-1)P, \quad s\in D^-,
$$
\beq
\Gs^+(s)=\left\{\fr{2}{\sqrt{\pi}K^+(s)}[X^+(s)]^{-1}X^+(-1)+I\right\}\fr{P}{s+1}, \quad s\in D^+,
\label{3.33}
\eeq
where $I$ is the unit matrix. Notice that when $\Ga=\fr{\pi}{2}$ the polynomial $f(s)$ vanishes, and we have
\beq
X^+(-1)=B^+(-1)
\left(\begin{array}{cc}
1+\Gb^+(-1)\cos\Ga & -2 \Gb^+(-1)\sin\Ga \\
0 & 1-\Gb^+(-1)\cos\Ga \\
\end{array}
\right).
\label{3.33'}
\eeq
We specify now 
formulas for the stress intensity factors. 
Instead of formula (\ref{3.21}) we have
\beq
\Gs^+(s)\sim\fr{2}{\sqrt{\pi}}(-s)^{-1/2}X_\infty^{-1}X^+(-1)P.
\label{3.34}
\eeq
Then, similar to the general case, we deduce expressions for the stress intensity factors. They are
\beq
\left(\begin{array}{c}
K_{I }\\
K_{II }\\
\end{array}
\right)=2\sqrt{\fr{2b}{\pi}}X_\infty^{-1}X^+(-1)P.
\label{3.35}
\eeq
It is helpful to rewrite this relation in the form
\beq
b^{-1/2}K_I=D_{11}P_1+D_{12}P_2, \quad b^{-1/2}K_{II}=D_{21}P_1+D_{22}P_2.
\label{3.35.1}
\eeq
In Table 1, we present the coefficients $D_{ij}$, $i,j=1,2$, which depend on the angle $\Ga$ only.

\vspace{5mm}

Table 1. The values of the coefficients $D_{ij}$, $i,j=1,2$, for some values of the angle $\Ga$ for constant loading.

\vspace{2mm}
\begin{tabular}{|c|c|c|c|c|}
\hline
$\Ga$     & $D_{11}$    &  $D_{12}$   & $D_{21}$       & $D_{22}$  \\
\hline
$\pi/8$    &   5.680018&  -0.643452 &     -2.941355 &   2.658874  \\       
\hline
$\pi/4$     &  2.791583  &    -0.355505   &  -0.937330   &   2.148703   \\   
\hline
$\pi/3$         &   2.225628 &  -0.252636     &  -0.542341    &  1.993084  \\       
\hline
$\pi/2$          &    1.776778 &   -0.121477   & -0.202058  &  1.813571   \\    
\hline
$2\pi/3$          &  1.636400 &   -0.0462352  &  -0.0661812 &   1.707510  \\   
\hline
$3\pi/4$         &  1.610895   & -0.0237185   &  -0.0322008 & 1.669096  \\    
\hline
$7\pi/8$         &   1.597337   & -0.00500464 &   -0.00639132    &     1.625163     \\    
\hline

\end{tabular}

\vspace{5mm}

When the constant loads applied to the crack faces are the same and constant, $\Gs_\Gt=-P_1$,  $\Gs_\Gt=-P_1$, $0<r<b$, $P_1=1$, the stress intensity 
factors tend to the same value as $\Ga\to \pi$ (Fig. 1). As $\Ga\to 0$, the magnitude of both factors is growing,   $K_I\to+\infty$, while  $K_{II}\to-\infty$.
Fig. 2 shows the stress intensity factors dependence on the angle $\Ga$ when $P_1=1$ and $P_2=0$. In this case $K_{II}\to 0$ as $\Ga\to \pi$. When $\Ga\to 0$,
the tendency of the factors $K_I$ and $K_{II}$ to go $+\infty$ and $-\infty$, respectively, preserves.

\vspace{2mm}

\subsection{Eigen-solutions}

As before, the wedge is $W=\{0<r<\infty, -\pi<\Gt<\Ga\}$,  the  crack is
$\{0<r<b, \Gt=0^\pm\}$, and  the wedge boundaries $\Gt=-\pi$ and $\Gt=\Ga$ are
free of traction. What is different now is that the crack faces are free of traction, while
at infinity, the stresses decay as $r^s$, $-1<s<0$. We call this problem $\CP^t$.
To apply the method described in  the previous sections, we
represent the solution to this problem, $\{\Gs_\Gt^t,\tau^t_{r\Gt}, \Gs^t_r\}$,
in the form
\beq
 \Gs^t_\Gt= \Gs_\Gt^\circ+\Gs_\Gt, \quad 
 \tau^t_{r\Gt}=\tau^\circ_{r\Gt}+\tau_{r\Gt},\quad
 \Gs^t_r= \Gs_r^\circ+\Gs_\Gt.
\label{3.36}
\eeq 
 Here, $\{\Gs_\Gt^\circ,\tau^\circ_{r\Gt}, \Gs^\circ_r\}$ is the solution to the 
 problem without a crack, while $\{\Gs_\Gt,\tau_{r\Gt}, \Gs_r\}$ is the solution 
 to the problem formulated in Section 2 with the stresses decaying at infinity as $r^{-1}$
 and with the boundary conditions on the crack faces given by
\beq
\Gs_\Gt(r,0^\pm)=-\Gs_\Gt^\circ(r,0), \quad \tau_{r\Gt}(r,0^\pm)=-\tau^\circ_{r\Gt}(r,0),  \quad 0<r<b.
\label{3.37}
\eeq 
We call the former problem $\CP^\circ$ and the latter problem $\CP$.

Since the total angle of the wedge is $\Ga+\pi\in(\pi,2\pi)$, the problem $\CP^\circ$
belongs to the class $N$
of singular plane problems [Cher]. Its solution can be found my the method of separation of variables. Assume that 
\beq
\Gs^\circ_\Gt(r,\Gt)= r^{s} k_\Gt(\Gt), \quad \tau^\circ_{r\Gt}(r,\Gt)= r^{s} k_{r\Gt}(\Gt),\quad
\Gs^\circ_r(r,\Gt)= r^{s} k_r(\Gt).
\label{3.38}
\eeq
On substituting the representations (\ref{3.38}) into the equilibrium equations (\ref{2.4}) and
the compatibility condition (\ref{2.5}) we reduce the system to the fourth order ordinary differential equation for the function $k_\Gt$
\beq 
k_\Gt^{IV}+[s^2+(s+2)^2]k''_\Gt+s^2(s+2)^2k_\Gt=0
\label{3.39}
\eeq
and the relations 
\beq
k_{r\Gt}=-\fr{k'_\Gt}{s+2}, \quad k_r=\fr{1}{s+1}\left(k_\Gt+\fr{k_\Gt''}{s+2}\right).
\label{3.39'}
\eeq
On solving the equation for $k_\Gt$ we find
$$
k_\Gt=C_1\cos s\Gt+C_2\cos(s+2)\Gt+C_3\sin s\Gt+C_4\sin(s+2)\Gt,
$$
\beq
k_{r\Gt}=s(s+2)^{-1}C_1\sin s\Gt+C_2\sin(s+2)\Gt-s(s+2)^{-1}C_3\cos s\Gt-C_4\cos(s+2)\Gt,
\label{3.40}
\eeq
where $C_1,\ldots,C_4$ are free constants. We next satisfy the homogeneous boundary conditions (\ref{2.1}) and find out that
if $\mu=s+1\in(0,1)$ is a root of the equation $d(s,\pi+\Ga)=0$ or, equivalently,
\beq
\sin^2\mu(\pi+\Ga)-\mu^2\sin^2\Ga=0,
\label{3.41}
\eeq
then one of the four constants $C_j$ say, $C_4$, is free and the others are expressed as
$$
C_1=\fr{C_4}{d_0(\mu)}[(\mu+1)\sin \mu\pi+\sin\mu(\pi+2\Ga)-\mu\sin(\mu\pi+2\Ga)],
$$
$$
C_2=\fr{C_4}{d_0(\mu)}[(\mu-1)\sin \mu\pi -\sin\mu(\pi+2\Ga)-\mu\sin(\mu\pi-2\Ga)],
$$
\beq
C_3=\fr{C_4}{d_0(\mu)}[(\mu+1)\cos \mu\pi-\cos\mu(\pi+2\Ga)-\mu\cos(\mu\pi+2\Ga)],
\label{3.42}
\eeq
where
\beq
d_0(\mu)=(\mu-1)\cos\mu\pi-\mu\cos(\mu\pi-2\Ga)+\cos\mu(\pi+2\Ga).
\label{3.43}
\eeq
On substituting these expressions into (\ref{3.40}) and using (\ref{3.39'}) and (\ref{3.38})
we determine the eigen-solutions $\{\Gs_\Gt^\circ,\tau^\circ_{r\Gt}, \Gs^\circ_r\}$ everywhere in the wedge.

Since the constant $C_4$ is free, we can rename the expression
\beq
k_\Gt^\circ=\fr{4C_4}{d_0(\mu)}\mu\sin\mu\pi\sin^2\Ga
\label{3.44}
\eeq
and write the boundary conditions (\ref{3.37}) as
\beq
\Gs_\Gt(r,0^\pm)=-k_\Gt^\circ r^{\mu-1}, \quad  
\tau_{r\Gt}(r,0^\pm)=-k_\Gt^\circ k_*r^{\mu-1}, \quad 0<r<b,
\label{3.45}
\eeq
where
$$
k_*=\fr{1}{4\mu(\mu+1)\sin\mu\pi\sin^2\Ga}[\mu(\mu+1)\cos(\mu\pi-2\Ga)
+\mu(\mu-1)\cos(\mu\pi+2\Ga)
$$
\beq
-2\cos\mu(\pi+2\Ga)-2(\mu^2-1)\cos\mu\pi].
\label{3.46}
\eeq
Having determined  the boundary conditions we can evaluate the Mellin 
transform of the load and obtain
\beq
g^-(s)=\fr{b^{\mu-1}P^*}{s+\mu}, \quad P^*= k_\Gt^\circ\left(\begin{array}{c}
1\\
k_*\\
\end{array}
\right).
\label{3.47}
\eeq
Similar to the case of the constant load considered in Section \ref{const}
the solution is found to be
$$
\Gc^-(s)=-\fr{4K^+(-\mu)b^{\mu-1}}{(s+\mu)K^-(s)}[X^-(s)]^{-1}X^+(-\mu)P^*, \quad s\in D^-,
$$
\beq
\Gs^+(s)=-\left\{\fr{K^+(-\mu)}{K^+(s)}[X^+(s)]^{-1}X^+(-\mu)-I\right\}\fr{P^*}{s+\mu}, \quad s\in D^+,
\label{3.48}
\eeq
The counterparts of the expressions (\ref{3.35}) in the case of the eigen-solutions 
have the form
\beq
\fr{1}{\sqrt{b}k_\Gt^\circ}\left(\begin{array}{c}
K_{I }\\
K_{II }\\
\end{array}
\right)=-\sqrt{2}K^+(-\mu)X_\infty^{-1}X^+(-\mu)\left(\begin{array}{c}
1\\
k_*\\
\end{array}
\right).
\label{3.49}
\eeq
Notice that the right hand-side in (\ref{3.49}) depends on the angle $\Ga$ only.

Analysis of the characteristic equation (\ref{3.41}) for the parameter $\mu$ shows that in the interval
$(0,1)$, if $0<\Ga<\Ga_*$, $\Ga_*\approx 0.431\pi$, there is only one root $\mu$ of equation (\ref{3.41}). Otherwise, if $\Ga_*<\Ga<\pi$, equation
(\ref{3.41}) has two roots, $\mu$ and $\mu_0$, in the interval $(0,1)$. Denote the smallest one $\mu$ and call the associated solution
 (\ref{3.38}) with $s=\mu-1$ the first eigen-solution.  We name the solution corresponding to the root $\mu_0$, $\mu<\mu_0<1$, the second eigen-solution.
 In comparison to the second solution the first solution has a stronger singularity at the edge and decays faster  at infinity.

We plot the stress singularity factors associated with the first and second eigen-solutions in Fig. 3 and 4, respectively. It is seen from Fig. 3 that in the former case 
both factors tend to infinity as $\Ga\to 0$. As $\Ga\to\pi$, $K_I\to 2.50663$ and $K_{II}\to 0$. We plot the factors $K_I$ and $K_{II}$ for the case of the second eigen-solution
for the values of $\Ga$ when the second solution exists. The factor $K_{II}$ is negative and its magnitude rapidly grows to infinity as $\Ga\to\pi$, while the factor $K_I$
is slowly growing as $\Ga$  grows. 

\begin{figure}[t]
\centerline{
\scalebox{0.6}{\includegraphics{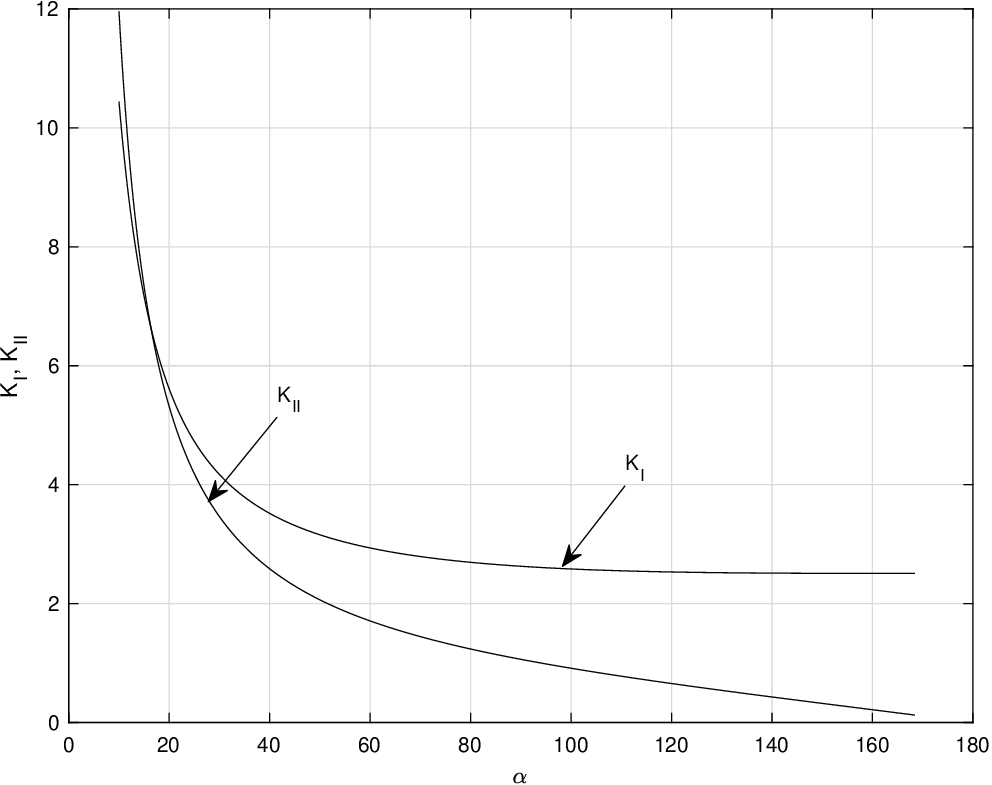}}}
\caption{
Stress intensity factors $K_I$ and $K_{II}$  versus the angle $\Ga$  when $b=1$, $k^\circ_\Gt=1$, and the loads correspond to the first eigen-solution: $p_\Gt(r)=r^{\mu-1}$, $p_{r\Gt}=k_* r^{\mu-1}$.}
\label{fig3}
\end{figure} 

\begin{figure}[t]
\centerline{
\scalebox{0.6}{\includegraphics{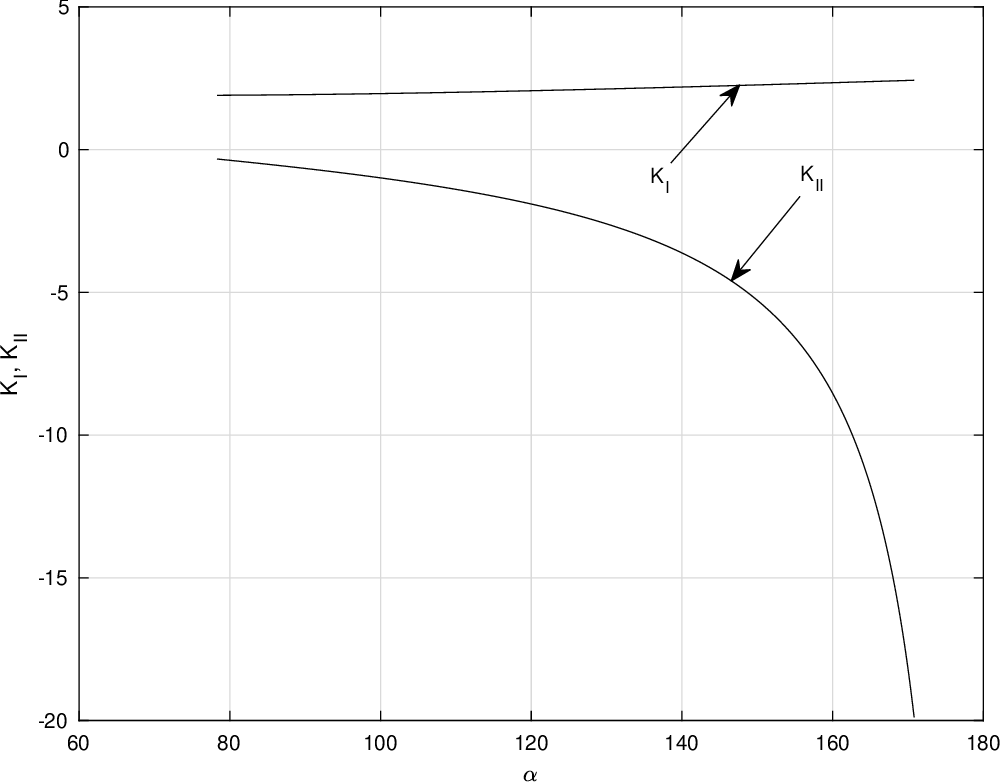}}}
\caption{
The stress intensity factors $K_I$ and $K_{II}$  versus the angle $\Ga$  when $b=1$, $k^\circ_\Gt=1$, and the loads correspond to the second eigen-solution are $p_\Gt(r)=r^{\mu_0-1}$, $p_{r\Gt}=k_* r^{\mu_0-1}$.}
\label{fig4}
\end{figure} 

Similar to the case of constant loading we represent formula (\ref{3.49}) as
\beq
\fr{K_I}{\sqrt{b}k_\Gt^\circ}=D_{11}+D_{12}k_*, \quad \fr{K_{II}}{\sqrt{b}k_\Gt^\circ}=D_{21}+D_{22}k_*,
\label{3.50}
\eeq
and give the values of $\mu$, $k_*$, and $D_{ij}$ for some values of the angle $\Ga$ corresponding to the first eigen-solution in Table 2.
In Table 3, we list the values of the second root $\mu_0$ in the interval $(0,1)$ and the stress intensity factors associated with the second eigen-solution for some
angles $\Ga$.

\vspace{5mm}

Table 2. The parameters $\mu$ and $k_*$ and the coefficients $D_{ij}$, $i,j=1,2$, for some values of the angle $\Ga$ in the case of the first eigen-solution.

\vspace{2mm}
\begin{tabular}{|c|c|c|c|c|c|c|}
\hline
$\Ga$     & $\mu$     & $k_*$ & $D_{11}$    &  $D_{12}$    & $D_{21}$     & $D_{22}$    \\
\hline
$\pi/8$   & 0.800766  &  2.721964 &  7.355120   &  -0.823447 &   -4.084327 &   3.219028  \\       
\hline
$\pi/4$   & 0.673583  & 1.340535 &    4.063677  &    -0.559928  &    -1.667266  &   2.956174  \\   
\hline
$\pi/3$   & 0.615731  & 0.9619090 &   3.363845 &  -0.445636  &   -1.094527  &    2.913322 \\       
\hline
$\pi/2$   &  0.544484 & 0.5430756 &   2.764929  &   -0.253395 &  -0.4847170   &    2.848868  \\    
\hline
$2\pi/3$  &  0.512221 &  0.302550 &  2.563696 &   -0.106728  &  -0.1738695  &  2.741192  \\   
\hline
$3\pi/4$  & 0.505010  &   0.213268 &   2.527096  &   -0.0565586 &    -0.0867114 &    2.676210 \\    
\hline
$7\pi/8$  &  0.500608 & 0.100301 &   2.508530   &     -0.0123458  &  -0.0175835  &  2.581612  \\    
\hline

\end{tabular}

\vspace{5mm}

Table 3. The parameter $\mu_0$ and the factors $\hat K_I=b^{-1/2}K_I/ k_\Gt^\circ$
and  $\hat K_{II}=b^{-1/2}K_{II}/ k_\Gt^\circ$  associated with the second eigen-solution
for some values of the angle $\Ga$.

\vspace{2mm}
\begin{tabular}{|c|c|c|c|}
\hline
$\Ga$     & $\mu_0$     & $\hat K_I$ & $\hat K_{II}$      \\
\hline
$0.435 \pi$   & 0.993452  &  1.904858 & -0.327423    \\       
\hline
$\pi/2$   &      0.908529  &   1.929324  &   -0.658165  \\    
\hline
$2\pi/3$  &    0.730901  &     2.064823 &    -1.904087  \\   
\hline
$3\pi/4$  &    0.659702 &     2.160225.  &   -3.055759  \\    
\hline
$7\pi/8$  &     0.570712  &    2.326037   &  -7.455055   \\    
\hline
$15\pi/16$  &   0.533251  &    2.415440  & -16.046723    \\    
\hline

\end{tabular}

\setcounter{equation}{0}

\section{Internal crack in a wedge along one of its sides}\label{int}

\subsection{Solution of the order-4 vector Riemann-Hilbert problem}
 
As in Section \ref{edge}, we take $\Ga_2=\pi$, while $\Ga_1=\Ga$ is an arbitrary angle in the range $(0,\pi)$. What is different now is that $a>0$ that is the crack is internal.
The problem is equivalent to the order-4 vector Riemann-Hilbert problem (\ref{2.16}) with the  matrix $G(s)$ given by (\ref{3.3}).
 We  split  the matrix $G(s)$ by means of the relation (\ref{3.5}),  factorize $\cot\pi s$ and the matrix $G_0(s)$ by formulas (\ref{3.7}) and (\ref{3.8}), respectively, and take into account that 
 the matrix $G(s)$ commutes with the Khrapkov factors of the matrix $G_0(s)$
\beq
G(t)X^-(t)=X^-(t)G(t), \quad G(t)X^+(t)=X^+(t)G(t),\quad t\in\GO.
\label{5.1}
\eeq
This enables us to rewrite the order-4 vector Riemann-Hilbert problem as follows:
$$
K^-(t)X^-(t)[\Gc^-(t)-\Gd^{t+1}G(t)\Gs^-(t)]=4K^+(t)X^+(t)[\Gs^+(t)+g^-(t)],
$$
\beq
\fr{[X^+(t)]^{-1}}{K^+(t)}[\Gc^+(t)-\Gd^{-t-1}G(t)\Gs^+(t)]=
\fr{4[X^-(t)]^{-1}}{K^-(t)}[\Gs^-(t)+g^+(t)], \quad t\in\GO.
\label{5.2}
\eeq 
For this model, we confine ourselves to considering  constant loading only, $p_\Gt(r)=P_1$, $p_{r\Gt}(r)=P_2$, $P_j=\const$.  Then
 the vectors $g^-(s)$ and $g^+(s)$
have the form
\beq
g^-(s)=\fr{\Gd^{s+1}-1}{s+1}P, \quad g^+(s)=\fr{1-\Gd^{-s-1}}{s+1}P,\quad 
P=\left(\begin{array}{c}P_1\\ P_2\\
\end{array}\right).
\label{5.3}
\eeq
Splitting the functions $g^\pm(t)$  into two parts, with $\Gd^{\pm t\pm 1}$ and without,
we transform the system (\ref{5.2}) as
$$
K^-(t)X^-(t)\left\{\Gc^-(t)-\Gd^{t+1}G(t)\left[\Gs^-(t)+\fr{P}{t+1}\right]\right\}
=4K^+(t)X^+(t)\left[\Gs^+(t)-\fr{P}{t+1}\right],
$$
\beq
\fr{[X^+(t)]^{-1}}{K^+(t)}\left\{\Gc^+(t)-\Gd^{-t-1}G(t)\left[\Gs^+(t)-\fr{P}{t+1} \right]\right\}=
\fr{4[X^-(t)]^{-1}}{K^-(t)}\left[\Gs^-(t)+\fr{P}{t+1}
\right], \quad t\in\GO.
\label{5.4}
\eeq 
First we remove the simple pole at $s=-1$ in the right hand-side of the first equation and 
the order-2 pole at $s=-1$ in the left hand-side
of the second equation. For this, we need the following asymptotic expansions
at $s=-1$:
$$
G(s)=\fr{4}{\pi(s+1)}G_0(-1)+O(s+1), \quad s\to -1,
$$
$$
[X^+(s)]^{-1}=[X^+(-1)]^{-1}+(s+1)\fr{d}{ds}[X^+(-1)]^{-1}+O((s+1)^2), \quad s\to -1,
$$$$
\fr{1}{K^+(s)}=-\fr{\sqrt{\pi}}{2}
[1-2(1-\log 2)(s+1)]+O((s+1)^2), \quad s\to -1.
$$
\beq
\Gd^{-s-1}=1-(s+1)\log\Gd+O((s+1)^2),  \quad s\to -1.
\label{5.5}
\eeq
Denote $s_1=1$. In addition to the simple poles at $\pm s_1\in D^\mp$ , the entries of the matrix $G(s)$ have a simple pole at $s_0=0\in D^-$ and simple poles at the 
points $\pm m\in D^\mp$  ($m=2,3,\ldots$) and at the complex-conjugate zeros  $\pm \Gs_m\in D^\mp$ of the function $d(s,\Ga)$. 
We denote all these zeros by $s_m$ and order them such that $\R s_m\le \R s_{m+1}$ ($m=2,3,\ldots$).
Aiming to eliminate these poles we introduce the vector-functions
\beq
\GL^+(s)=\sum_{m=0}^\infty \fr{A_m^+}{s-s_m}, \quad   
\GL^-(s)=\sum_{m=1}^\infty \fr{A_m^-}{s+s_m}, \quad 
A_m^\pm=\left(
\begin{array}{c}
A_{1m}^\pm\\
A_{2m}^\pm\\
\end{array}
\right).
\label{5.6}
\eeq
The series $\GL^\pm(s)$ are holomorphic vector-functions in the half-planes $D^\pm$ and have simple poles
of the matrix $G(s)$ in the domains $D^\mp$. The residues of these vector-functions are to 
be found. It will be shown later that the series  $\GL^\pm(s)$ converge uniformly and absolutely in the 
domains 
\beq
D^\pm_\Gve=C\setminus\cup_{m=m_\mp}^\infty D^\mp_\Gve(\pm s_m), 
\label{5.7}
\eeq
where $C$ is the complex $s$-plane, $m_-=0$, $m_+=1$,  and $D_\Gve^\mp(\pm s_m)$ are discs of a small radius $\Gve>0$
centered at the points $\pm s_m\in D^\mp$. On removing the poles mentioned, using the asymptotics
of the factors $K^\pm(s)$ at infinity
\beq
K^\pm(s)\sim\mp(\mp s)^{-1/2}, \quad s\to\infty, \quad s\in D^\pm,
\label{5.8} 
\eeq
and since $\Gc^\pm(s)=O(s^{-1/2})$, $\Gs^\pm(s)=O(s^{-1/2})$. $s\to\infty$, $s\in D^\pm$,
we arrive at the following relations which continue each other analytically into the whole complex plane:
$$
K^-(s)X^-(s)\left\{\Gc^-(s)-\Gd^{s+1}G(s)\left[\Gs^-(s)+\fr{P}{s+1}\right]\right\}-\fr{N_1^-}{s+1}-\GL^+(s)
$$$$
=4K^+(s)X^+(s)\left[\Gs^+(s)-\fr{P}{s+1}\right]-\fr{N_1^-}{s+1}-\GL^+(s)=0,\quad s\in C,
$$
$$
\fr{[X^+(s)]^{-1}}{K^+(s)}\left\{\Gc^+(s)-\Gd^{-s-1}G(s)\left[\Gs^+(s)-\fr{P}{s+1} \right]\right\}
-\fr{N_1^+}{s+1}-\fr{N_2^+}{(s+1)^2}-\GL^-(s)
$$
\beq
=\fr{4[X^-(s)]^{-1}}{K^-(s)}\left[\Gs^-(s)+\fr{P}{s+1}
\right]-\fr{N_1^+}{s+1}-\fr{N_2^+}{(s+1)^2}-\GL^-(s)=C^\circ, \quad s\in C.
\label{5.9}
\eeq 
Here, $C^\circ=(C_1^\circ,C_2^\circ)$ is an arbitrary constant vector, while  $N_1^+$,
 $N_2^+$ and $N_1^-$ are constant vectors given by
$$
N_1^+=\fr{2}{\sqrt{\pi}}\left\{[X^+(-1)]^{-1}[2(1-\log 2)+\log \Gd]-\fr{d}{ds}[X^+(-1)]^{-1}\right\}G_0(-1)P,
$$
\beq
N_2^+=-\fr{2}{\sqrt{\pi}}G_0(-1)[X^+(-1)]^{-1}P,\quad 
N_1^-=\fr{8}{\sqrt{\pi}}X^+(-1)P.
\label{5.10}
\eeq
The system (\ref{5.9}) enables us to write  representation formulas for the solution of the vector Riemann-Hilbert problem (\ref{2.16}). They read
$$
\Gs^+(s)=\fr{P}{s+1}+\fr{1}{4K^+(s)}[X^+(s)]^{-1}\GS_1(s), \quad s\in D^+,
$$
$$
\Gs^-(s)=-\fr{P}{s+1}+\fr{K^-(s)}{4}X^-(s)\GS_2(s), \quad s\in D^-,
$$
$$
\Gc^-(s)=\fr{\Gd^{s+1}}{4}G(s)K^-(s)X^-(s)\GS_2(s)+\fr{1}{K^-(s)}[X^-(s)]^{-1}\GS_1(s), \quad s\in D^-,
$$
\beq
\Gc^+(s)=\fr{\Gd^{-s-1}}{4 K^+(s)}G(s)[X^+(s)]^{-1}\GS_1(s)
+K^+(s)X^+(s)\GS_2(s), \quad s\in D^+,
\label{5.11}
\eeq
where the vector-functions $\GS_1(s)$ and $\GS_2(s)$ are
\beq
\GS_1(s)=\fr{N_1^-}{s+1}+\GL^+(s), 
\quad 
\GS_2(s)=C^\circ +\fr{N_1^+}{s+1}+\fr{N_2^+}{(s+1)^2}+\GL^-(s).
\label{5.12}
\eeq
It is directly verified that $\Gc^-(s)=\Gd^{s+1}\Gc^+(s)$, and the vector-functions $\Gc^\pm(s)$
satisfy the first equation in (\ref{2.16}). The vector-functions $\Gs^\pm(s)$ are holomorphic
in the half-planes $D^\pm$. The simple poles of the vector-functions $\Gc^\pm(s)$ at the poles of the
entries of the matrix $G(s)$ become removable singularities if and only if the
vectors $A^\pm_m$ solve the following infinite system of  linear algebraic equations:
$$
A^+_n=\Gd^{s_n+1}\GD^-_n\left[C^\circ +\fr{N_1^+}{s_n+1}+\fr{N_2^+}{(s_n+1)^2}
+\sum_{m=1}^\infty\fr{A_m^-}{s_n+s_m}\right], \quad n=0,1,\ldots,
$$
\beq
A^-_n=\Gd^{s_n-1}\GD^+_n\left(h_nN_1^-
-\sum_{m=0}^\infty\fr{A_m^+}{s_n+s_m}\right), \quad n=1,2,\ldots,
\label{5.13}
\eeq
where
$$
\GD_n^+=-\fr{1}{4[K^+(-s_n)]^2}\mathop{\mathrm{res}}_{s = -s_n}G(s)\{[X^+(-s_n)]^{-1}\}^2, \quad n=1,2,\ldots,
$$
$$
\GD_n^-=-\fr{1}{4}[K^-(s_n)]^2\mathop{\mathrm{res}}_{s = s_n}G(s)[X^-(s_n)]^2,\quad n=0,1,\ldots,
$$
$$
h_1=-2(1-\log 2)[G_0(-1)]^{-1}+X^+(-1)\fr{d}{ds}[X^+(-1)]^{-1}, \quad h_n=\fr{1}{1-s_n}, \quad n\ge 2,
$$$$
\mathop{\mathrm{res}}_{s =\pm \Gs_n}G(s)
=\fr{1}{2\Gs_n\sin^2\Ga-\Ga\sin 2\Ga \Gs_n}
$$
$$
\times\left(
\begin{array}{cc}
-\sin 2\Ga \Gs_n-\Gs_n\sin 2\Ga & -2\Gs_n(\pm \Gs_n-1)\sin^2\Ga\\
2\Gs_n(\pm \Gs_n+1)\sin^2\Ga & -\sin 2\Ga \Gs_n+\Gs_n\sin 2\Ga\\
\end{array}
\right), 
$$
\beq
\mathop{\mathrm{res}}_{s =\pm n}G(s)
=\fr{4}{\pi}G_0(\pm n), \quad 
n=0,1,2,\ldots.
\label{5.14}
\eeq
It is clear that the coefficients $A_n^\pm=O(\Gd^{n})$, $n\to \infty$, and not only 
 the series (\ref{5.6}) converge
uniformly and absolutely with an exponential rate but also an approximate solution of the infinite system (\ref{5.13})
converges exponentially to its exact solution.

\subsection{Definition of the vector $C^\circ$}
 
 The solution has to satisfy the additional conditions (\ref{2.3'}) or, equivalently,
 the vector $\Gc^-(s)$ has to meet the condition
 \beq
 \Gc^-(0)=0.
 \label{5.15}
 \eeq
 This vector equation is employed to fix the unknown vector $C^\circ=(C_1^\circ, C_2^\circ)$.
 To satisfy the condition (\ref{5.15}), we write asymptotic expansions 
for small $s$
$$
\Gd^{s+1}=\Gd[1+s\log\Gd +O(s^2)],
$$
$$
[K^-(s)]^{\pm 1}=\pi^{\pm 1/2}(1\mp 2s\log 2)+O(s^2), 
$$
$$
\GL^-(s)=a_1^- -s a_2^- +O(s^2), \quad \GL^+(s)=\fr{A_0^+}{s}-a_1^+ +O(s),
$$
$$
G(s)=\fr{4}{\pi s}G_0(0)+O(s),
$$
\beq
[X^-(s)]^{\pm 1}=[X^-(0)]^{\pm 1}+s\fr{d}{ds}[X^-(0)]^{\pm 1}+O(s^2), \quad s\to 0.
\label{5.16}
\eeq
Here,
\beq
a_j^\pm=\sum_{n=1}^\infty \fr{1}{s_n^j}A_n^\pm, \quad j=1,2.
\label{5.17}
\eeq
On substituting these expressions into the representation formula for the vector-function 
$\Gc^-(s)$ in (\ref{5.11}) we obtain
$$
\Gc^-(s)=\fr{1}{s\sqrt{\pi}}\{[X^-(0)]^{-1}A_0^++\Gd G_0(0)X^-(0)(C^\circ+a_1^-+N_1^++N_2^+)\}
+\fr{\Gd}{\sqrt{\pi}}G_0(0)
$$$$
\times\left\{-X^-(0)(a_2^-+N_1^++2N_2^+)+\left[\fr{d}{ds}X^-(0)+\log\fr{\Gd}{4}X^-(0)\right]
(C^\circ+a_1^-+N_1^++N_2^+)\right\}
$$
\beq
+\fr{1}{\sqrt{\pi}}\left(2\log 2 [X^-(0)]^{-1}+\fr{d}{ds}[X^-(0)]^{-1}\right)A_0^+
+\fr{1}{\sqrt{\pi}}[X^-(0)]^{-1}(N_1^--a_1^+)+O(s), \quad s\to 0.
\label{5.18}
\eeq
It follows from the first equation in (\ref{5.13}) for $n=0$
that the factor of $1/s$ in (\ref{5.18}) is equal to 0. Therefore the condition (\ref{5.15})
is met if and only if the vector $C^\circ$ solve the vector equation
$$
\left[\fr{d}{ds}X^-(0)+\log\fr{\Gd}{4}X^-(0)\right]
(C^\circ+a_1^-+N_1^++N_2^+)=X^-(0)(a_2^-+N_1^++2N_2^+)
$$
\beq
-\fr{1}{\Gd} [G_0(0)]^{-1}\left\{\left(2\log 2 [X^-(0)]^{-1}+\fr{d}{ds}[X^-(0)]^{-1}\right)A_0^+
+[X^-(0)]^{-1}(N_1^--a_1^+)\right\}.
\label{5.19}
\eeq

Determine next the principal term of the asymptotic expansion of the vector $C$
as $\Gd\to 0$. We need this term to derive the asymptotics of the stress intensity
factors when the crack is close to the wedge vertex that is when $a\to 0$.
Analysis of the system (\ref{5.13}) and formula (\ref{5.19}) for small $\Gd$ shows that
\beq
\Gd\log\Gd G_0(0) X^-(0) C^\circ\sim -[X^-(0)]^{-1} N_1^-, \quad \Gd\to 0.
\label{5.20}
\eeq
Employing  formula (\ref{5.10}) for $N_1^-$ and since
 $G_0(0)=X^+(0)[X^-(0)]^{-1}$,  we have
\beq
X^-(0)X^+(0)C^\circ\sim-\fr{8}{\sqrt{\pi}\Gd\log\Gd}X^+(-1)P, \quad \Gd\to 0.
\label{5.21}
\eeq
Prove next that 
the product
$X^-(0)X^+(0)$ is the unit matrix.
Notice that $f(0)=1$ and the functions $\GD(t)$ and $f^{-1/2}(t)\Ge(t)$ are even.
By the Sokhotski-Plemelj forrmulas, we deduce
\beq
B^\pm(0)=e^{\pm\fr14\log\GD(0)}, \quad \Gb^\pm(0)=\pm\fr12\Ge(0),
\label{5.22}
\eeq
and therefore
\beq
c_+(0^\pm)=c_0\pm s_0\cos\Ga, 
\quad c_-(0^\pm)=c_0\mp s_0\cos\Ga, 
\quad
s_+(0^\pm)=s_-(0^\pm)=\mp s_0\sin\Ga.
\label{5.23}
\eeq
where
\beq
c_0=\cosh\fr12\Ge(0), \quad s_0=\sinh\fr12\Ge(0).
\label{5.24}
\eeq
The matrices $X^\pm(0)$ have the form
\beq
X^\pm(0)=e^{\pm\fr14\log\GD(0)}
\left(
\begin{array}{cc}
c_0\pm s_0\cos\Ga & \mp s_0\sin\Ga\\
\mp s_0\sin\Ga & c_0\mp s_0\cos\Ga \\
\end{array}
\right)
\label{5.25}
\eeq
and their product is the unit matrix. The asymptotics of  the vector $C^\circ$ is simplified,
\beq
C^\circ\sim-\fr{8}{\sqrt{\pi}\Gd\log\Gd}X^+(-1)P, \quad \Gd\to 0.
\label{5.26}
\eeq

\subsection{Stress intensity factors and the potential energy}

Introduce the stress intensity factors in the standard. way
$$
\Gs_\Gt(r,0)\sim\fr{K_I^-}{\sqrt{2\pi(a-r)}}, \quad 
\tau_{r\Gt}(r,0)\sim\fr{K_{II}^-}{\sqrt{2\pi(a-r)}}, \quad r\to a-0,
$$
\beq
\Gs_\Gt(r,0)\sim\fr{K_I^+}{\sqrt{2\pi(r-b)}}, \quad 
\tau_{r\Gt}(r,0)\sim\fr{K_{II}^+}{\sqrt{2\pi(r-b)}}, \quad r\to b+0.
\label{5.27}
\eeq 
On applying the abelian theorem to the integrals $\Gs^-(s)$ and $\Gs^+(s)$ given by (\ref{2.15})
we obtain their asymptotics at infinity
$$
\Gs^-(s)\sim \fr{s^{-1/2}}{\sqrt{2a}}\left(\begin{array}{c}
K_I^-\\
K_{II}^-\\
\end{array}
\right), \quad s\to\infty, \quad s\in D^-,
$$
\beq
\Gs^+(s)\sim \fr{(-s)^{-1/2}}{\sqrt{2b}}\left(\begin{array}{c}
K_I^+\\
K_{II}^+\\
\end{array}
\right), \quad s\to\infty, \quad s\in D^+.
\label{5.28}
\eeq
On the other side, the asymptotics of the vector-functions $\Gs^\pm(s)$ is defined from the analysis
of the solution (\ref{5.11})
$$
\Gs^-(s)\sim \fr14 s^{-1/2} X_\infty C^\circ, \quad s\to\infty, \quad s\in D^-,
$$
\beq
\Gs^+(s)\sim \fr14 (-s)^{-1/2} [X_\infty]^{-1}\left(N_1^-+\sum_{m=0}^\infty A_m^+\right),
 \quad s\to\infty, \quad s\in D^+,
 \label{5.29}
 \eeq
where 
\beq
X_\infty=
\left(\begin{array}{cc}
\cos q & -\sin q\\
\sin q & \cos q\\
\end{array}
\right)
\label{5.32}
\eeq
and $q$ is given by (\ref{3.28}).
By comparing formulas (\ref{5.28}) and (\ref{5.29}) it is straightforward to find the stress intensity factors
\beq
\left(\begin{array}{c}
K_I^-\\
K_{II}^-\\
\end{array}
\right)=\fr12\sqrt{\fr{a}{2}}X_\infty C^\circ,\quad 
\left(\begin{array}{c}
K_I^+\\
K_{II}^+\\
\end{array}
\right)=\fr12\sqrt{\fr{b}{2}}[X_\infty]^{-1}\left(N_1^-+\sum_{m=0}^\infty A_m^+\right).
\label{5.33}
\eeq

It is of interest to find the asymptotics of the  stress intensity factors when $b$ is fixed, while
$a\to 0$ or, equivalently, when $\Gd\to 0$. Using the asymptotics (\ref{5.26}) of the vector $C^\circ$
and since $A_m^+=O(1/\log\Gd)$, $\Gd\to 0$, we deduce
$$
\left(\begin{array}{c}
K_I^+\\
K_{II}^+\\
\end{array}
\right)\sim 2\sqrt{\fr{2b}{\pi}}[X_\infty]^{-1}X^+(-1)P, \quad \Gd\to 0,
$$
\beq
\left(\begin{array}{c}
K_I^-\\
K_{II}^-\\
\end{array}
\right)\sim -\fr{2\sqrt{2b}}{\sqrt{\pi\Gd}\log\Gd}
X_\infty X^+(-1)P, \quad \Gd\to 0.
\label{5.34}
\eeq

\begin{figure}[t]
\centerline{
\scalebox{0.6}{\includegraphics{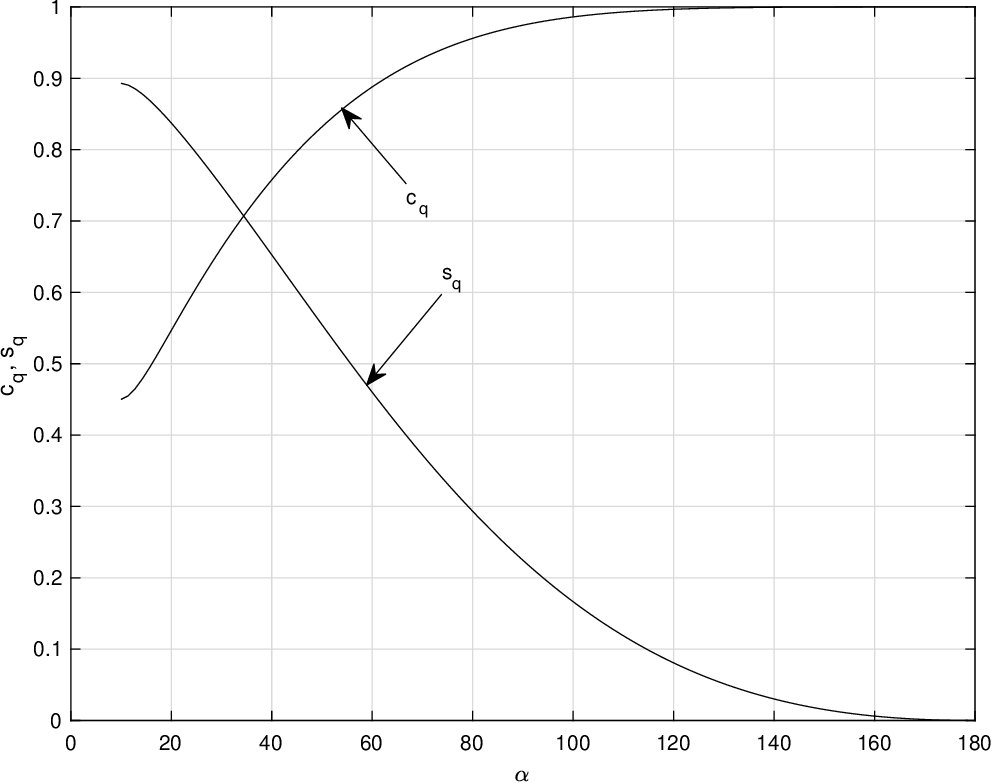}}}
\caption{
Entries $c_q=\cos 2q$ and  $c_s=\sin 2q$  of the matrix $Q$ versus the angle $\Ga$.}
\label{fig5}
\end{figure} 

Denote by $Q$ the matrix $X^2_\infty$ that is
\beq
Q=\left(\begin{array}{cc}
\cos 2q & -\sin 2q\\
\sin 2q & \cos 2q\\
\end{array}
\right)
\label{5.35}
\eeq
and rewrite the asymptotic formulas (\ref{5.34}) as
\beq
\left(\begin{array}{c}
K_I^+\\
K_{II}^+\\
\end{array}
\right)\sim
\left(\begin{array}{c}
K_I\\
K_{II}\\
\end{array}
\right),
\quad 
\left(\begin{array}{c}
K_I^-\\
K_{II}^-\\
\end{array}
\right)\sim -\fr{1}{\sqrt{\Gd}\log\Gd} Q\left(\begin{array}{c}
K_I\\
K_{II}\\
\end{array}
\right),
 \quad \Gd\to 0,
\label{5.36}
\eeq
where 
\beq
\left(\begin{array}{c}
K_I\\
K_{II}\\
\end{array}
\right)=2\sqrt{\fr{2b}{\pi}}[X_\infty]^{-1}X^+(-1)P
\label{5.37}
\eeq
are the stress intensity factors
at the tip $r=b$ when the crack is at the vertex of the wedge that is when $a=0$,
$P=(P_1,P_2)$, and  the loads 
$\Gs_\Gt(r,0)=-P_1$, $\tau_{r\Gt}(r,0)=-P_2$, $0<r<b$,
are constant. 

\begin{figure}[t]
\centerline{
\scalebox{0.6}{\includegraphics{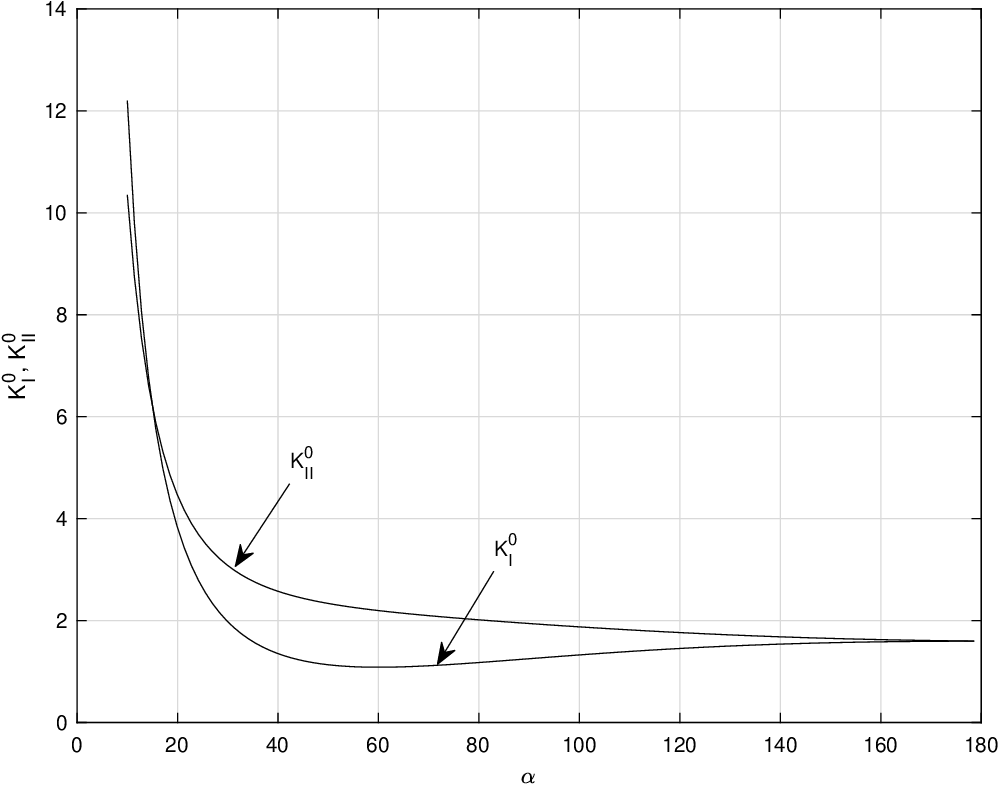}}}
\caption{
Case $a\to 0$: factors $K_I^0$ and $K_{II}^0$  versus the angle $\Ga$  when $b=1$ and  $P_1=P_2=1$.}
\label{fig6}
\end{figure}

The variations of the entries $c_q=\cos 2q$ and  $s_q=\cos 2q$ of the matrix $Q$ with the angle $\Ga$ are shown in Fig. 5. 
The factors 
\beq
K_I^0=\cos 2q K_I-\sin 2q K_{II}, \quad K_{II}^0=\sin 2q K_I+\cos  2q K_{II},
\label{5.38}
\eeq
in the asymptotic relations
\beq
K^-_I\sim -\fr{K_I^0}{\sqrt{\Gd}\log\Gd}, \quad K^-_{II}\sim -\fr{K_{II}^0}{\sqrt{\Gd}\log\Gd},  \quad \Gd\to 0,
\label{5.39}
\eeq
as functions of the angle $\Ga$ and when $b=1$ and  $P_1=P_2=1$,
are plotted in Fig. 6. 

In the final part of this section we determine the potential energy $\Gd U$ released when the crack $r\in(a,b)$ extends first to $r\in (a,b+h)$ and then
when the original crack $r\in(a,b)$ extends to the left, $r\in(a-h,b)$, $h>0$. In the former case, we have
\beq
\Gd U^+ \sim\fr12\int_b^{b+h}
\{\Gs_\Gt(r,0)\Gd [u_\Gt](r,0)+\tau_{r\Gt}\Gd [u_r](r,0)\}dr,
\label{5.40}
\eeq
where $[u_\Gt]+\Gd[u_\Gt]$ and $[u_r]+\Gd[u_r]$ are the displacement jumps from the upper side to the lower side.
We have
\beq
\Gc(r)\sim \fr{L^+}{\sqrt{b-r}}, \quad r\to b-0,
\label{5.41}
\eeq
where $L^+$ is a vector with constant components. By integrating we deduce
\beq
[u](r,0)\sim \fr{2}{E}L^+\sqrt{b-r},  \quad r\to b-0,
\label{5.42}
\eeq
where $[u](r,0)=([u_\Gt], [u_{r\Gt}])(r,0)$. Analysis of representation (\ref{5.11}) of the function $\Gc^-(s)$ as $s\to \infty$,
the use of the abelian theorem and formula (\ref{5.33}) yield 
\beq
L^+=\fr{4}{\sqrt{2\pi}}\left(\begin{array}{c}
K^+_I\\
K^+_{II}\\
\end{array}
\right),
 \quad 
[u](r,0)\sim \fr{4}{E}\sqrt{\fr{2}{\pi}}\left(\begin{array}{c}
K^+_I\\
K^+_{II}\\
\end{array}
\right)\sqrt{b-r}, \quad r\to b-0.
\label{5.43}
\eeq
Therefore, when the crack extends from $r=b$ to $r=b+h$, we have
\beq
\Gd[u](r,0)\sim \fr{4}{E}\sqrt{\fr{2}{\pi}}\left(\begin{array}{c}
K^+_I\\
K^+_{II}\\
\end{array}
\right)\sqrt{b+h-r}, \quad r\to h+b-0.
\label{5.44}
\eeq
On substituting the expressions (\ref{5.27}) and (\ref{5.44}) into formulas (\ref{5.40}), evaluating the integral  and repeating the analysis for the crack tip $r=a$ we find
\beq
\Gd U^\pm \sim\fr{h}{E}[(K_I^\pm)^2+(K_{II}^\pm)^2], \quad h\to 0,
\label{5.45}
\eeq
where $\Gd U^-$ is the potential energy released when the crack extends from 
 $r=a$ to $r=a-h$, while $b$ is fixed.
When the crack tip $r=a$ is close to the wedge vertex, that is when $\Gd\to 0$, we employ the asymptotic relations 
(\ref{5.36}) and deduce
\beq
\Gd U^- \sim\fr{h}{E\Gd\log^2\Gd} (K_I^2+K_{II}^2), \quad \Gd U^+ \sim\fr{h}{E} (K_I^2+K_{II}^2),
\quad h\to 0,  \quad \Gd\to 0.
\label{5.46}
\eeq
Here, as before, $K_I$ and $K_{II}$ are the stress intensity factors at the crack tip $r=b$ when  $a=0$.

\setcounter{equation}{0}

\section{General case}

In the general case of a crack $\{a<r<b, \Gt=0^\pm\}$ in a wedge
 $W=\{-\Ga_2<\Gt<\Ga_1, 0<r<\infty\}$, the entries $G_{jm}$ of the matrix $G(s)$
are given by (\ref{2.12}), and the matrix does not admits an exact factorization by the methods currently
available. The aim of this section is to show that even in the general case under the assumption that a factorization exists
and when $a\to 0$ ($\Gd\to 0$),  the stress intensity factors 
at the crack tip $r=a$ asymptotically behave as the corresponding factors at the tip $r=b$ multiplied by $-\Gd^{-1/2}/\log\Gd$ and a matrix that depends 
on the angles $\Ga_1$ and $\Ga_2$  and is independent of the loading.  First we split the matrix $G(s)$ 
into the dominant scalar 
that yields the singularity of the solution at infinity and the ``regular" matrix 
\beq
G(s)=4\cot\pi s G^\circ(s).
\label{6.1}
\eeq
The diagonal entries of the matrix $G^\circ(s)$ are positive as $s\to 0$ and
tend to 1 as $s=\Gs+i\tau$ and $\tau\to\pm\infty$, while the off diagonal elements
decay exponentially as $\tau\to\pm\infty$. The determinant
of the matrix $G^\circ(s)$ tends to 1 as   $\tau\to\pm\infty$ and is positive
as $s\to 0$. It is directly verified that 
the increment of the argument of the determinant of the matrix $G^\circ(s)$ as $s$ traverses the contour $\GO$
is zero, that is the total index of the matrix $G^\circ(s)$ is equal to 0.
We also assume that there exist the left and right standard factorizations \cite{goh} of the matrix $G(t)$
\beq
G(t)=4 K^+(t)X_L^+(t)[K^-(t)X_L^-(t)]^{-1}=4[K^-(t)X_R^-(t)]^{-1}K^+(t)X_R^+(t),
\quad t\in\GO,
\label{6.2}
\eeq
with the left and right partial indices being zeros. Here,  $K^\pm(s)$ are the functions given by (\ref{3.7}), and
the right factors $X_R^\pm(t)$ of the matrix $G^\circ(t)$ are expressed through the left factors $X_R^\pm(t)$ of the inverse 
matrix $[G^\circ(t)]^{-1}$ by the formula
\beq
X_R^\pm(t)=[\hat X_L^\pm(t)]^{-1}, \quad t\in\GO.
\label{6.2'}
\eeq
 
On following the procedure of Section \ref{int}
we express the solution in terms of the matrices $X_R^\pm(s)$ and $X_L^\pm(s)$
$$
\Gs^+(s)=\fr{P}{s+1}+\fr{1}{4K^+(s)}[X_R^+(s)]^{-1}\GS_1(s), \quad s\in D^+,
$$
$$
\Gs^-(s)=-\fr{P}{s+1}+\fr{K^-(s)}{4}X_L^-(s)\GS_2(s), \quad s\in D^-,
$$
$$
\Gc^-(s)=\fr{\Gd^{s+1}}{4}G(s)K^-(s)X_L^-(s)\GS_2(s)+\fr{1}{K^-(s)}[X_R^-(s)]^{-1}\GS_1(s), \quad s\in D^-,
$$
\beq
\Gc^+(s)=\fr{\Gd^{-s-1}}{4K^+(s)}G(s)[X_R^+(s)]^{-1}\GS_1(s)
+K^+(s)X_L^+(s)\GS_2(s), \quad s\in D^+,
\label{6.3}
\eeq
where the vector-functions $\GS_1(s)$ and $\GS_2(s)$ are given in (\ref{5.12}) and 
the constant vector $C^\circ$ solves the equation 
$$
\left[\fr{d}{ds}X_L^-(0)+\log\fr{\Gd}{4}X_L^-(0)\right]
(C^\circ+a_1^-+N_1^++N_2^+)=X_L^-(0)(a_2^-+N_1^++2N_2^+)
$$
\beq
-\fr{1}{\Gd} [G_0(0)]^{-1}\left\{\left(2\log 2 [X_R^-(0)]^{-1}+\fr{d}{ds}[X_R^-(0)]^{-1}\right)A_0^+
+[X_R^-(0)]^{-1}(N_1^--a_1^+)\right\}.
\label{6.4}
\eeq
Similar to the case $\Ga_2=\pi$ we derive the asymptotics of the vector $C^\circ$ as $\Gd\to 0$.
We have
\beq
C\sim -\fr{1}{\Gd\log\Gd}[X_L^+(0)]^{-1}[X_R^-(0)]^{-1} N_1^-, \quad \Gd\to 0.
\label{6.5}
\eeq
This result enables us to derive the counterpart of the asymptotic relations (\ref{5.36}). 
We have
\beq
\left(\begin{array}{c}
K_I^+\\
K_{II}^+\\
\end{array}
\right)\sim
\left(\begin{array}{c}
K_I\\
K_{II}\\
\end{array}
\right)=2\sqrt{\fr{2b}{\pi}}[X_{R\infty}]^{-1}X_R^+(-1)P
\label{6.6}
\eeq
and
\beq
\left(\begin{array}{c}
K_I^-\\
K_{II}^-\\
\end{array}
\right)\sim -\fr{1}{\sqrt{\Gd}\log\Gd} Q\left(\begin{array}{c}
K_I\\
K_{II}\\
\end{array}
\right),
 \quad \Gd\to 0,
\label{6.7}
\eeq
where the matrix-factor  $Q$ is 
\beq
Q=X_{L\infty}[X_L^+(0)]^{-1}[X_R^-(0)]^{-1}X_{R\infty},
\label{6.8}
\eeq
and the matrices $X_{R\infty}$ and $X_{L\infty}$ are the limits as $s\to\infty$ of the matrices $X_R(s)$ and 
$X_L(s)$,  respectively, while $K_I$ and $K_{II}$ are the stress intensity factors at the tip $b$ when $a=0$.

\setcounter{equation}{0}

\section{Internal crack orthogonal to the boundary of a half-plane}

Assume now that $\Ga_1=\Ga_2=\pi/2$, $0<a<b$, and the crack faces are subjected to normal loading, $\Gs_\Gt=-p_\Gt(r)$, $a<r<b$.
The order-4 Riemann-Hilbert problem (\ref{2.16}) is simplified and becomes an order-2 problem
$$
\Gc^+(s)=\Gd^{-s-1}\Gc^-(s),
$$
\beq
\Gs^+(s)=L(s)\Gc^-(s)-\Gd^{s+1}\Gs^-(s)-g^-(s),\quad s\in\GO,
\label{7.1}
\eeq
 where the following notations are accepted:
 $$
 L(s)=\fr{\sin^2\fr{\pi s}{2}-s^2}{2\sin\pi s}, 
$$$$
\Gc^-(s)=\int_\Gd^1\Gc_1(br)r^s dr, \quad \Gc^+(s)=\int_1^{1/\Gd}\Gc_1(ar)r^s dr, 
$$$$
\Gs^-(s)=\int_0^1\Gs_\Gt(ar,0)r^s dr, \quad \Gs^+(s)=\int_1^{\infty}\Gs_\Gt(br,0)r^s dr.
$$
\beq
g^-(s)=-\int_\Gd^1 p_\Gt(br)r^s dr, \quad 
g^+(s)=-\int_1^{1/\Gd}p_\Gt(ar)r^s dr,
\label{7.2}
\eeq
and $\Gc_1(r)$ is the jump of the tangential derivative of the normal displacement introduced in (\ref{2.3}). Note that the problem can equivalently be formulated by means of the singular integral 
equation
\beq
\fr{1}{4\pi}\int_a^b\Gc_1(\Gr)
\left(\fr{1}{\Gr-r}+\fr{r^2+4r\Gr-\Gr^2}{(\Gr+r)^3}\right)d\Gr=p_\Gt(r), \quad a<r<b,
\label{7.3}
\eeq
and reduced to the vector Riemann-Hilbert problem (\ref{7.1}) by employing the Mellin transform and the Mellin convolution theorem.
Following the scheme of Section \ref{int} we factorize the function $L(s)$,
\beq
L(s)=\fr{\tan\pi s}{4}L_0(s), \quad L_0(s)=1-\fr{s^2}{\sin^2\fr{\pi s}{2}},
\label{7.4}
\eeq
in the form
$$
L(t)=-\fr14\fr{L^+(t)}{L^-(t)}, \quad t\in\GO,
$$$$
L^\pm(s)=a^\pm(s)X^\pm(s),\quad 
a^+(s)=\fr{\GG(\fr12-\fr{s}{2})}{\GG(-\fr{s}{2})}, \quad a^-(s)=\fr{\GG(1+\fr{s}{2})}{\GG(\fr12+\fr{s}{2})},
$$
$$
X^\pm(s)=\exp\left\{\fr{1}{2\pi i}\int_\GO\fr{\log L_0(t)dt}{t-s}\right\}
$$
\beq=\exp\left\{-\fr{s}{\pi}\int_0^\infty\log\left( 1-\fr{\tau^2}{\sinh^2\fr{\pi \tau}{2}}\right)\fr{d\tau}{\tau^2+s^2}\right\}, \quad s\in D^\pm.
\label{7.5}
\eeq
In the case of constant  loading, $\Gs_\Gt(r,0)=-P=\const$, $a<r<b$, after substituting the factorization formulas into (\ref{7.1}) we deduce
$$
-\fr{\Gc^-(s)}{4L^-(s)}
+\fr{\Gd^{s+1}}{4L^-(s)L(s)}\left[\Gs^-(s)+\fr{P}{s+1}\right]+\fr{p_+}{s+1}-\GL^+(s)
$$$$
=\fr{1}{L^+(s)}\left[\Gs^+(s)-\fr{P}{s+1}\right]
+\fr{p_+}{s+1}-\GL^+(s)=0,\quad s\in C,
$$$$
L^-(s)\left[\Gs^-(s)+\fr{P}{s+1}\right]-\fr{p_-}{s+1}-\GL^-(s)
$$
\beq
=-\fr14 L^+(s)\Gc^+(s)+\fr{\Gd^{-s-1}L^+(s)}{4L(s)}\left[\Gs^+(s)-\fr{P}{s+1}\right]-\fr{p_-}{s+1}-\GL^-(s)=C^\circ, \quad s\in C.
\label{7.6}
\eeq
Here, 
$$
p_+=\fr{P}{L^+(-1)}, \quad p_-=\fr{\pi P L^+(-1)}{4},
$$
\beq
\GL^+(s)=\fr{A_0^+}{s}+\sum_{m=1}^\infty \fr{A_m^+}{s-s_m}, \quad   
\GL^-(s)=\sum_{m=1}^\infty \fr{A_m^-}{s+s_m},
\label{7.7}
\eeq
and $\pm s_m\in D^\mp$ ($m=1,2,\ldots$) are complex-conjugate zeros of the function $s^2-\sin^2 \fr{\pi s}{2}$. The coefficients  $A_0^+$, $A_m^\pm$ ($m=1,2,\ldots$),
and the constant $C^\circ$ are to be determined.

Relations (\ref{7.6}) enable us to express the solution of the vector Riemann-Hilbert problem in terms of the series (\ref{7.7})
$$
\Gs^+(s)=\fr{P}{s+1}+L^+(s)\left[\GL^+(s)-\fr{p_+}{s+1}\right],
$$$$
\Gs^-(s)=-\fr{P}{s+1}+\fr{1}{L^-(s)}\left[\GL^-(s)+C^\circ+\fr{p_-}{s+1}\right],
$$
$$
\Gc^-(s)=\fr{\Gd^{s+1}}{L^-(s)L(s)} [\GL^-(s)+C^\circ +\fr{p_-}{s+1}]
-4L^-(s)\left[\GL^+(s)-\fr{p_+}{s+1}\right],
$$
\beq
\Gc^+(s)=\fr{\Gd^{-s-1}L^+(s)}{L(s)} \left[\GL^+(s)-\fr{p_+}{s+1}\right]
-\fr{4}{L^+(s)}\left[\GL^-(s)+C^\circ +\fr{p_-}{s+1}\right].
\label{7.8}
\eeq
It is immediately verified that $\Gc^+(s)=\Gd^{-s-1}\Gc^-(s)$.
In general, for arbitrary chosen coefficients  $A_0^+$ and $A_m^\pm$ ($m=1,2,\ldots$), the functions
$\Gc^\pm(s)$ have inadmissible simple poles in the half-planes $D^\pm$, respectively. The poles become
removable singularities if the coefficients solve the following system of linear algebraic equations:
$$
A_n^+=\Gd^{s_n+1}\GD^-_n\left(C^\circ+\fr{p_-}{s_n+1}+\sum_{m=1}^\infty\fr{A_m^-}{s_n+s_m}\right),\quad n=0,1,\ldots,
$$
\beq
A_n^-=\Gd^{s_n-1}\GD^+_n\left(\fr{p_+}{s_n-1}-\sum_{m=0}^\infty\fr{A_m^+}{s_n+s_m}\right),\quad n=1,2,\ldots,
\label{7.9}
\eeq
where
\beq
\GD^\pm_n=\fr{[L^\pm(\mp s_n)]^{\pm 2}\sin\pi s_n}{\pi\sin\pi s_n-4s_n}, \quad n=1,2,\ldots,\quad \GD_0^-=\fr{2\pi}{(\pi^2-4)[L^-(0)]^2}, \quad s_0=0.
\label{7.10}
\eeq
We shift the contour $\GO$ to the imaginary axis and use the Sokhotski-Plemelj formula for the limit values of the function $X(s)$ on the contour $\GO_0$
\beq
X^\pm(t)=\exp\left\{\pm\fr12\log L_0(t)+\fr{1}{2\pi i} P.V.\int_{\GO_0}\log L_0(\tau)\fr{d\tau}{\tau-t}\right\}, \quad t\in\GO_0,
\label{7.11}
\eeq
where $\GO_0=\{\R t=0,  -\infty<\I t<\infty\}$.
Passing to the limit $t\to 0^+$ we find 
\beq
X^-(0)=\fr{\pi}{\sqrt{\pi^2-4}}, \quad L^-(0)=\sqrt{\fr{\pi}{\pi^2-4}}, \quad \GD^-_0=2.
\label{7.12}
\eeq
The constant $C^\circ$ is fixed by the condition $\Gc^-(0)=0$. Similar to Section \ref{int} we employ the first equation 
for $n=0$ in (\ref{7.9}) and derive the following equation for the constant $C^\circ$:
\beq
2\Gd[(C^\circ +a_1^-+p_-)(\log\Gd -L_0)-a_2^--p_-]+a_1^++p_+-L_0A_0^+=0,
\label{7.13}
\eeq
where
\beq
L_0=\fr{1}{L^-(0)}\fr{d}{ds}L^-(0), \quad a_j^\pm=\sum_{m=1}^\infty\fr{A_m^\pm}{s_m^j}.
\label{7.14}
\eeq
To compute $L_0$, we integrate by parts in the integral in formula (\ref{7.11}) and differentiate $X^-(t)$ 
\beq
\fr{d}{dt}[X^-(t)]=X^-(t)\left[-\fr{L'_0(t)}{2L_0(t)}+\fr{1}{2\pi i} P.V.\int_{\GO_0}\fr{L'_0(\tau)}{L_0(t)}\fr{d\tau}{\tau-t}\right].
\quad  t\in\GO_0.
\label{7.15}
\eeq
By  passing to the limit $t\to 0$ and taking into account that
\beq
L_0'(\tau)=\fr{2\tau}{\sin^2\fr{\pi\tau}{2}}\left(-1+\fr{\pi\tau}{2}\cot\fr{\pi\tau}{2}\right)\sim -\fr23\tau, \quad \tau\to 0,
\label{7.16}
\eeq
we derive 
\beq
\fr{d}{ds}X^-(0)=X^-(0)X_1,
\label{7.17}
\eeq
where
\beq
X_1=\fr{2}{\pi}\int_0^\infty\fr{\fr{\pi\tau}{2}\coth\fr{\pi\tau}{2}-1}{\tau^2-\sinh^2\fr{\pi\tau}{2}}d\tau.
\label{7.18}
\eeq
The integrand is bounded and has the limit $\fr13\pi^2(4-\pi^2)^{-1}$ at the point 0 and exponentially decays at infinity.
Now, since 
\beq
\fr{d}{ds}a^-(s)=\fr{\GG(1+\fr{s}{2})}{2\GG(\fr12+\fr{s}{2})}
\left[\psi\left(1+\fr{s}{2}\right)-\psi\left(\fr12+\fr{s}{2}\right)\right]\sim\fr{\log 2}{\sqrt{\pi}},\quad s\to 0,
\label{7.19}
\eeq
we have the following formula for $L_0$:
\beq
L_0=X_1+\log 2.
\label{7.20}
\eeq

It is convenient to eliminate the constant $C^\circ$ from the infinite system (\ref{7.9}). We represent the unknown coefficients in the form
\beq
A^\pm_m=C A_{0m}+A_{1m}^\pm
\label{7.21}
\eeq
and denote
\beq
f_{0n}^-=1, \quad f_{1n}^-=\fr{p_-}{s_n+1}, \quad f_{0n}^+=0, \quad f_{1n}^+=\fr{p_+}{s_n-1}.
\label{7.22}
\eeq
Then the new coefficients $A_{0m}^\pm$ and  $A_{1m}^\pm$ are solutions of two separate systems of linear algebraic equations
$$
A_{jn}^+=\Gd^{s_n+1}\GD^-_n\left(f^-_{jn}+\sum_{m=1}^\infty\fr{A_{jm}^-}{s_n+s_m}\right),\quad n=0,1,\ldots,
$$
\beq
A_{jn}^-=
\Gd^{s_n-1}\GD^+_n
\left(f_{jn}^+ -\sum_{m=0}^\infty\fr{A_{jm}^+}{s_n+s_m}
\right),
\quad n=1,2,\ldots.
\label{7.23}
\eeq
Due to the factors $\Gd^{s_n+1}$ and $\Gd^{s_n-1}$ the rate of convergence of an approximate solution to the exact one is exponential.

On substituting the representations (\ref{7.21}) into the equation for the constant $C^\circ$ (\ref{7.13}) and denoting
\beq
a_{jk}^\pm=\sum_{m=1}^\infty \fr{A_{jm}^\pm}{s_m^k}
\label{7.24}
\eeq
we express the constant $C^\circ$ through the solutions of the systems (\ref{7.23}) 
\beq
C^\circ=
\fr{a_{11}^+ +p_+-  L_0A_{10}^+-2\Gd[(\log\Gd-L_0)(a_{11}^-+p_-)-a_{12}^- -p_-]}
{2\Gd[(\log\Gd -L_0)(a_{01}^- +1)-a_{02}^-]-a_{01}^+ +L_0A_{00}^+}.
\label{7.25}
\eeq
 
 Finally, we determine the stress intensity factors
 \beq
 K^-_I=\lim_{r\to a^-}\sqrt{2\pi(a-r)}\Gs_\Gt(r,0), \quad   K^+_I=\lim_{r\to b^+}\sqrt{2\pi(r-b)}\Gs_\Gt(r,0).
\label{7.26}
\eeq
From one side, by the abelian transform for the Mellin integrals,
\beq
\Gs^-(s)\sim \fr{s^{-1/2}}{\sqrt{2a}}K_I^-, \quad s\in D^-, 
\quad \Gs^+(s)\sim
\fr{(-s)^{-1/2}}{\sqrt{2b}}K_I^+, \quad s\in D^+, \quad s\to\infty.
\label{7.27}
\eeq
On the other side, from the representation formulas (\ref{7.8}), the asymptotics of these functions at infinity have the form
$$
\Gs^-(s)\sim\left(\fr{s}{2}\right)^{-1/2}C,  \quad s\in D^-, \quad s\to\infty.
$$
\beq
\Gs^+(s)\sim \fr{(-s)^{-1/2}}{\sqrt{2}}
\left(p_+ -\sum_{m=1}^\infty A_m^+\right), \quad s\in D^+, \quad s\to\infty.
\label{7.28}
\eeq
Formulas (\ref{7.27}) and (\ref{7.28}), when combined, give the stress intensity factors at the crack tips $r=a$ and $r=b$ in the case of constant loading
\beq
K_I^-=2\sqrt{a}C, \quad K_I^+=\sqrt{b}\left(\fr{P}{L^+(-1)} -\sum_{m=0}^\infty A_m^+\right).
\label{7.29}
\eeq
It is of interest to find the asymptotics of the factors when $b$ is fixed and  $a\to 0$ or, equivalently, when $\Gd\to 0$. 
Analysis of equation (\ref{7.13}) shows that
\beq
C\sim -\fr{P}{2 L^+(-1)}\fr{1}{\Gd\log\Gd}, \quad \Gd\to 0.
\label{7.30}
\eeq
From the system (\ref{7.9}), 
\beq
\sum_{m=0}^\infty A_m^+=O\left(\fr{1}{\log\Gd}\right), \quad \Gd\to 0.
\label{7.31}
\eeq
Therefore 
\beq
K_I^+\sim K_I, \quad  K_I^-\sim -\fr{K_I}{\sqrt{\Gd}\log\Gd}, \quad \Gd\to 0,
\label{7.32}
\eeq
where $K_I$ is the stress intensity factor at the crack tip $r=b$ when $a=0$
\beq
K_I=\sqrt{\pi b}P\Gg, \quad \Gg=\exp\left\{-\fr{1}{\pi}\int_0^\infty \log\left(1-\fr{\tau^2}{\sinh^2\fr{\pi\tau}{2}}\right)\fr{d\tau}{\tau^2+1}\right\}\approx 1.1215222.
\label{7.33}
\eeq
This result is consistent with the value $\Gg=1.1215$ obtained in \cite{koi2}.

Finally, as in Section \ref{int}, we determine the potential energy $\Gd U^-$ released when the crack extends from 
 $r=a$ to $r=a-h$, while $b$ is fixed, and the corresponding value $\Gd U^+$ associated with the crack tip $r=b$,
 \beq
 \Gd U^\pm\sim\fr{h}{E}(K_I^\pm)^2, \quad h\to 0.
 \label{7.34}
 \eeq
When the crack is close to the boundary of the half-plane, we have
\beq
 \Gd U^-\sim\fr{h}{E\Gd\log\Gd}K_I^2,  \quad \Gd U^+\sim\fr{h}{E}K_I^2, \quad h\to 0, \quad \Gd\to 0.
 \label{7.35}
 \eeq

\setcounter{equation}{0}
  
\section{Conclusion}

We have analyzed model problems of a crack in a wedge when the crack $a<r<b$ is an edge slit ($a=0$) or an internal cut ($a>0$). 
Our attention has been focussed on the case when the crack $\{0\le a<r<b, \Gt=0\}$ is located in the wedge $\{0<r<\infty, -\pi<\Gt<\Ga\}$,
$0<\Ga<\pi$.
The solution \cite{khr2} to this problem of an edge crack $a=0$  was obtained by matrix factorization of the kernel of the 
order-2 vector    Wiener-Hopf problem. The results \cite{khr2} also include approximate formulas for the weight matrix neededed
for the stress intensity factors. We have modified the factorization \cite{khr2} by first splitting the kernel
into a dominant scalar function responsible for the solution singularities and a ``regular" matrix whose factorization 
is more convenient for numerical purposes. We have considered two cases of the loading and determined the stress intensity factors
by quadratures. The first case is standard and concerns
constant normal and tangential loads applied to the crack faces. The second one uses the eigen-solution of the 
singular problem of elasticity of a wedge whose angle $\Ga+\pi\in (\pi,2\pi)$ \cite{che}. The wedge boundary and the crack faces are assumed to be free
of traction, while at infinity, the stresses decay as $r^{\mu-1}$, where $\mu\in(0,1)$ is a real root of the associated 
characteristic equation of the wedge without a crack.  

In the case of an internal crack, we reduced the model to an order-4 vector Riemann-Hilbert problem. For its solution,
we have proposed a method that uses the $2\times 2$ matrix factorization arising in the model of the edge crack and requiring
solving an infinite system of linear algebraic equations of the second kind. What is important is that the rate of convergence of the approximate
solution to the exact one is exponential. 
We have derived the stress intensity factors $K_{I,II}^+$ and $K_{I,II}^-$ at the crack tips $r=b$ and $r=a$, respectively.
It has been discovered that
 \beq
\left(\begin{array}{c}
K_I^-\\
K_{II}^-\\
\end{array}
\right)\sim -\fr{1}{\sqrt{\Gd}\log\Gd} Q\left(\begin{array}{c}
K_I\\
K_{II}\\
\end{array}
\right),
 \quad \Gd=\fr{a}{b}\to 0,
\label{8.1}
\eeq
where $Q$ is a matrix whose entries are independent of the loading and functions of the angle $\Ga$ only and $K_{I,II}=K_{I,II}^+|_{a=0}$.
We have established that, when the crack $\{a<r<b, \Gt=0\}$ advances to the left, $r=a-h$ ($h>0$ and small) and $a$ is close to the wedge
 vertex, then the potential energy increment $\Gd U$ has the form
 \beq
\Gd U \sim\fr{h}{E\Gd\log^2\Gd} (K_I^2+K_{II}^2),
\quad h\to 0,  \quad \Gd\to 0,
\label{8.2}
\eeq  
 where $E$ is the Young modulus. 
 
 In the general case of the problem, when the crack $\{a<r<b, \Gt=0\}$ is in an arbitrary wedge $\{0<r<\infty, -\Ga_1<\Gt<\Ga_2\}$,
the associated $2\times 2$ matrix kernel cannot be factorized by the methods available in the literature. However, if the left and right
factorizations become available, then the asymptotic formula (\ref{8.1}) is also valid with the matrix $Q$ expressible 
through not explicitly defined Wiener-Hopf factors and independent of the loading. Finally, we note that in the case of an internal crack orthogonal to the boundary 
of a half-plane and close to the boundary, formula (\ref{8.1}) is simplified, and the stress intensity factors at the tips $a$
and  $b$, $K_I^-$ and $K_I|_{a=0}$, respectively, are connected by
\beq
K_I^-\sim -\fr{K_I}{\sqrt{\Gd}\log\Gd}, \quad \Gd\to 0.
\label{8.3}
\eeq

\end{document}